\documentclass[pdflatex,sn-mathphys-num]{sn-jnl}% Math and Physical Sciences Numbered Reference Style
\usepackage{graphicx}%
\usepackage{multirow}%
\usepackage{amsmath,amssymb,amsfonts}%
\usepackage{amsthm}%
\usepackage{mathrsfs}%
\usepackage[title]{appendix}%
\usepackage{xcolor}%
\usepackage{textcomp}%
\usepackage{manyfoot}%
\usepackage{booktabs}%
\usepackage{algorithm}%
\usepackage{algorithmicx}%
\usepackage{algpseudocode}%
\usepackage{listings}%
\usepackage{color,xcolor,ucs}
\usepackage{mathtools}   \usepackage{tikz} 
\usepackage{ amssymb }
\usepackage{extarrows} 
\usepackage{pgf,tikz}
\usepackage{float}
\usetikzlibrary{positioning}
\usetikzlibrary{shapes.geometric}
\usetikzlibrary{shapes.misc}
\usetikzlibrary{arrows}
\usepackage{caption}
\usepackage{mathrsfs}
\usetikzlibrary{arrows,shapes,automata,backgrounds,petri,positioning}
\usetikzlibrary{decorations.pathmorphing}
\usetikzlibrary{decorations.shapes}
\usetikzlibrary{decorations.text}
\usetikzlibrary{decorations.fractals}
\usetikzlibrary{decorations.footprints}
\usetikzlibrary{shadows}
\usetikzlibrary{calc}
\usetikzlibrary{spy}
\usepackage{amsmath}
\usepackage{array}
\usepackage{ amssymb }
\usepackage{braket}
\usepackage{qcircuit}
\usepackage{soul}
\usepackage{braket} 
\usepackage{relsize}

\usepackage{amsmath}
\usepackage{ amssymb }
\usepackage{braket}
\usepackage{qcircuit}
\usepackage{soul}
\usepackage{braket}
\theoremstyle{thmstyleone}%
\theoremstyle{thmstyletwo}%

\theoremstyle{thmstylethree}%

\begin{document}

\title[Article Title]{Phase transition of the long range Ising model in lower dimensions, for $d < \alpha \leq d + 1$, with a Peierls’ argument}


\author[1]{\fnm{Pete} \sur{Rigas}}\email{rigas.pete@gmail.com}

\affil[1]{\orgaddress{ \city{Newport Beach}, \postcode{92625}, \state{CA}, \country{United States}}}


\abstract{We extend previous results due to Ding and Zhuang in order to prove that a phase transition occurs for the long range Ising model in lower dimensions. By making use of a recent argument due to Affonso, Bissacot and Maia from 2022 which establishes that a phase transition occurs for the long range, random-field Ising model, from a suggestion of the authors we demonstrate that a phase transition also occurs for the long range Ising model, from a set of appropriately defined contours for the long range system, and a Peierls’ argument.}




\maketitle

\section{Introduction}

\subsection{Overview}

\noindent The random-field Ising model, RFIM, is a model of interest in statistical mechanics, not only for connections with the celebrated Ising model, through the phenomena of ferromagnetism [8], but also for connections with the random-field, long-range Ising model which was shown to exhibit a phase transition [1], correlation length lower bounds with the greedy lattice animal [3], a confirmation of the same scaling holding for the correlation length of the random-field Potts model [8], long range order [4], Monte Carlo studies [9], community structure [11], supersymmetry [13], and the computation of ground states [14]. To extend previous methods for proving that a phase transition occurs in the random-field, long-range Ising model besides only one region of $\alpha$ parameters dependent on the dimension $d$ of the lattice, we implement the argument for analyzing contours, provided in [1], for the contours provided in [2]. In comparison to arguments for proving that the phase transition occurs in [1], in which a variant of the classical Peierls’ argument is implemented by reversing the direction of the spins contained within contours $\gamma$, the contours described in [2] can be of use for proving that the phase transition for the random-field, long-range Ising model occurs for another range of $\alpha$ parameters, in which $d < \alpha \leq d+1$.

Beginning in the next section, after having defined the model, as well as connections that it shares with the random-field, and long-range Ising model, we introduce contour systems for the long range, random-field, and long range Ising models, from which we conclude with a Peierls’ argument for proving that a phase transition occurs.

\subsection{This paper's contributions}

\noindent The following work seeks to present a review of a recent argument, due to Affonso, Bissacot and Maia, which exhibited that a phase transition occurs for RFIM using a Peierls’ argument. Arguments of this form were originally put forth in seminal work due to Peierls’, as a means of correcting a conclusion, due to Ising, that the Ising model does not undergo a phase transition. For the RFIM, given a restriction on the parameter $\alpha$ which appears in the specific form of the coupling constants of the Hamiltonian, one can formulate that the RFIM undergoes a phase transition in higher dimensions, namely for $\alpha \geq d+1$. The authors of this effort, Affonso, Bissacot and Maia, suggest that their methods could remain applicable for demonstrating that the long range, RFIM undergoes a phase transition in lower dimensions, specifically for $\alpha$ satisfying $d < \alpha \leq d+1$. In such a lower dimensional setting, in comparison to the long range RFIM in higher dimensions, the phase transition states that there exists two probability measures which are distinct under $+$ and $-$ boundary conditions. Conveniently, the lower-dimensional RFIM has boundary conditions which can also be encoded through $+$ and $-$ spins as the higher-dimensional RFIM. Furthermore, under this assignment of boundary conditions, the methodology developed by Affonso, Bissacot and Maia, is applicable, through arguments consisting of the following steps: (1) formulating lower-dimensional, long range, coupling constants; (2) formulating a spin flip procedure, which as a lower-dimensional analog to such a procedure in higher dimensions, reverses the orientation of spins contained within the interior of suitable contours; (3) demonstrating that the probability of sampling a length, from the suitable contour system, with a length exceeding some threshold is exponentially small; (4) upper bounding the probability, under long-range interactions, of $\big\{ \sigma_0 \equiv -1 \big\}$, namely the spin at the origin being of a $-$ sign, with two exponentials scaling with respect to parameters taken to be sufficiently small. With respect to the long, or short, range parameter $\alpha$ which determines the typical scale over which a $+$, or $-$, spin influences the orientation of nearest neighbors, the RFIM having a phase transition regardless of the regime of parameters for $\alpha$ is enlightening, and unique to phase transitions of other models in the field, such as the FK-Ising model, which fails to have a continuous phase transition when $q>4$.

\subsection{Long range, random-field Ising model objects}

\noindent After having motivated the physical characteristics of phase transitions for models in Statistical Physics, we introduce objects for encoding boundary conditions on long-range probability measures. Despite the fact that boundary conditions are encoded identically for the lower-dimensional, and higher-dimensional, RFIM models, the range over which interactions between neighboring sites are expected to influence interactions throughout other portions of the lattice determines how finite volume measures are expected to behave asymptotically in weak finite volume. In particular, the lower-dimensional RFIM, as does the higher-dimensional RFIM, still experiences a phase transition regardless of the exponent introduced for defining the coupling constants.

To introduce the probability measure for the long range, random-field Ising model, first consider, for a finite volume $\Lambda \subsetneq \textbf{Z}^d$, with $\big| \Lambda \big| < + \infty$,

\begin{align*}
  \mathcal{H}^{\mathrm{LR}, \eta}_{\Lambda} \big( \sigma \big)  = - \underset{x,y \in \Lambda}{\sum} J_{x,y} \sigma_x \sigma_y - \underset{y \in \Lambda^c}{\underset{x \in \Lambda}{\sum}} J_{x,y} \sigma_x \eta_y \text{, }  
\end{align*}

\noindent corresponding to the Hamiltonian for the long-range Ising model, in which the spins in the first, and second summation, have coupling constants $\big\{ J_{xy} \big\}_{x,y \in \textbf{Z}^d}$, spins $\sigma_x$ and $\sigma_y$ in $\Lambda$, spin $\eta_y$ in $\Lambda^c$ for the boundary conditions, each of which is drawn from the spin-sample space $\Omega \equiv \big\{ -1 , 1 \big\}^{\textbf{Z}^d}$, with coupling constants,

\begin{align*}
   J_{xy} \equiv J \big| x - y \big|^{-\alpha} \text{, }  x \neq y \text{, }   
\end{align*}

\noindent for some strictly positive $J$, $\alpha > d$, and $J_{xy} = 0$ otherwise. The couplings for the Hamiltonian, in both the long-range, and random-field case introduced below, are also intended to satisfy,

\begin{align*}
\underset{|x|>1}{\underset{x \in \textbf{Z}^d}{\sum}} \big| x_i \big| J_{0,x} < J_{0,e_i}              \text{, }  
\end{align*}

\noindent in which the couplings are translation invariant, for every $1 \leq i \leq d$. In the presence of disorder through an external field, specifically through the iid family of Gaussian variables $\big\{ h_x \big\}_{x \in \textbf{Z}^d}$, the long-range, random-field Ising model Hamiltonian takes the form,

\begin{align*}
 \mathcal{H}^{\mathrm{LR-RF}, \eta}_{\Lambda}  \equiv  \mathcal{H}^{\mathrm{LR-RF}, \eta}_{\Lambda} \big( \sigma , h  \big) = \mathcal{H}^{\mathrm{LR}, \eta}_{\Lambda} \big( \sigma) - \underset{x \in \Lambda}{\sum}  \epsilon h_x \sigma_x \text{, }  
\end{align*}

\noindent which is also taken under boundary conditions $\eta$, for some strictly positive $\epsilon$. In the summation over $x \in \Lambda$ above besides the long-range Hamiltonian terms, the external field takes the form,

\[
h_x  \equiv  \text{ } 
\left\{\!\begin{array}{ll@{}>{{}}l}     h^{*}   &  \text{, if } 
x = 0  \text{, }  \\
h^{*} \big| x \big|^{-\delta}  &   \textit{, if } x \neq 0  \text{ , }  \\
\end{array}\right.
\]

\noindent for $\delta, h^{*}>0$. The corresponding Gibbs measure, 

\begin{align*}
 \textbf{P}_{\Lambda,\beta}  \big( \sigma , h  \big) 
\equiv \textbf{P}^{\mathrm{LR-RF},\eta}_{\Lambda,\beta }  \big( \sigma , h \big) \equiv  \frac{\mathrm{exp} \big( \beta \mathcal{H}^{\mathrm{LR-RF}, \eta}_{\Lambda}\big) }{Z^{\mathrm{LR-RF},\eta}_{\Lambda ,\beta} \big( h \big)  }    \text{, }  
\end{align*}

\noindent at inverse temperature $\beta >0$, has the partition function as the normalizing constant so that $\textbf{P} \big( \cdot \big)$ is a probability measure, with,

\begin{align*}
Z^{\mathrm{LR-RF},\eta}_{\Lambda ,\beta }  \big( h \big)  \equiv   \underset{x \in \Omega^{\eta}_{\Lambda} }{\sum}\mathrm{exp} \big(  \beta \mathcal{H}^{\mathrm{LR-RF}, \eta}_{\Lambda}\big) 
   \text{, }  
\end{align*}

\noindent over the sample space $\Omega^{\eta}_{\Lambda}$ of spins with boundary condition $\eta$ over $\Lambda$. Similarly, for the long range Ising model,

\begin{align*}
  \textbf{P}_{\Lambda,\beta}  \big( \sigma , h  \big) 
 \equiv \textbf{P}^{\eta}_{\Lambda} \big( \sigma , h \big) \equiv \textbf{P}^{\mathrm{LR},\eta}_{\Lambda,\beta }  \big( \sigma , h \big) \equiv  \frac{\mathrm{exp} \big(  \beta \mathcal{H}^{\mathrm{LR}, \eta}_{\Lambda}\big) }{Z^{\mathrm{LR},\eta}_{\Lambda ,\beta} \big( h \big)  } \equiv  \frac{\mathrm{exp} \big(  \beta \mathcal{H}_{\Lambda}\big) }{Z^{\eta}_{\Lambda ,\beta} \big( h \big) }  \text{, }  
\end{align*}

\noindent with,

\begin{align*}
Z^{\mathrm{LR},\eta}_{\Lambda ,\beta }  \big( h \big)  \equiv Z^{\eta}_{\Lambda ,\beta }  \big( h \big)  \equiv   \underset{x \in \Omega^{\eta}_{\Lambda} }{\sum}\mathrm{exp} \big(  \beta \mathcal{H}^{\mathrm{LR}, \eta}_{\Lambda}\big) 
   \text{. } 
\end{align*}

\noindent Equipped with $\textbf{P} \big( \cdot \big)$, the joint probability measure for the pair $\big( \sigma , h \big)$ is,

\begin{align*}
  \textbf{Q}^{\mathrm{LR-RF},{\eta}}_{\Lambda , \beta} \big( \sigma \in A , h \in B \big) \equiv \underset{B}{\int}            \textbf{P}^{\mathrm{LR-RF},\eta}_{\Lambda_,\beta } \big(  A \big)   \text{ }  \mathrm{d} \textbf{P}^{\mathrm{LR-RF}}_{\Lambda , \beta}  \big( h \big) \text{, }  \\ \textbf{Q}^{\mathrm{LR},{\eta}}_{\Lambda , \beta} \big( \sigma \in A , h \in B \big) \equiv \underset{B}{\int}            \textbf{P}^{\mathrm{LR},\eta}_{\Lambda_,\beta } \big(  A \big)   \text{ }  \mathrm{d} \textbf{P}^{\mathrm{LR},\eta}_{\Lambda , \beta}  \big( h \big) \equiv  \underset{B}{\int}            \textbf{P}^{\eta}_{\Lambda_,\beta } \big(  A \big)   \text{ }  \mathrm{d} \textbf{P}^{\eta}_{\Lambda , \beta}  \big( h \big)   \text{, }  
\end{align*}

\noindent under boundary conditions $\eta$, for measurable $A \subsetneq \Omega$ and $B \subsetneq \textbf{R}^{\Lambda}$ borelian, with density,

\begin{align*}
  \mathcal{D}^{\mathrm{LR-RF},\eta}_{\Lambda , \beta} \big( \sigma , h \big)   = \underset{u \in \Lambda}{\prod} \frac{1}{\sqrt{2 \pi}} \text{ }  \mathrm{exp} \big(  - \frac{h^2_u }{2}      \big) \text{ } \textbf{P}^{\mathrm{LR-RF}, \eta}_{\Lambda_,\beta } \big(  \sigma , h \big)  \text{, }     \\ \mathcal{D}^{\mathrm{LR},\pm}_{\Lambda , \beta} \big( \sigma , \eta \big)  \equiv  \mathcal{D}^{\pm}_{\Lambda , \beta} \big( \sigma , \eta \big)  = \underset{u \in \Lambda}{\prod} \frac{1}{\sqrt{2 \pi}} \text{ }  \mathrm{exp} \big(  - \frac{\eta^2_u }{2}      \big) \text{ } \textbf{P}^{\eta}_{\Lambda_,\beta } \big(  \sigma , \eta \big)  \text{, }   
\end{align*}

\noindent under $+$ boundary conditions. As a sequence of finite volumes $\Lambda_n$, with $\Lambda_n \subsetneq \Lambda$ and $ \big| \Lambda_n \big| < + \infty$, tends to $\textbf{Z}^d$ via a weak limit,

\begin{align*}
    \textbf{P}^{\mathrm{LR-RF},\eta}_{\beta } \big[  \omega    \big] =      \underset{n \longrightarrow + \infty }{\mathrm{lim}}  \textbf{P}^{\mathrm{LR-RF},\eta}_{\Lambda_n,\beta } \big[  \omega   \big] \text{, }  \\   \textbf{P}^{\mathrm{LR},\eta}_{\beta } \big[  \omega    \big] =      \underset{n \longrightarrow + \infty }{\mathrm{lim}}  \textbf{P}^{\mathrm{LR},\eta}_{\Lambda_n,\beta } \big[  \omega   \big] \text{, }  
\end{align*}

\noindent for a random-field, long-range Ising configuration $\omega \in \Omega^{\eta}_{\Lambda}$. From seminal work in [3], the authors of [1] extend work for proving that the phase transition for the random-field Ising model occurs, introduced in [3], surrounding a \textit{Peierls’ type argument} for demonstrating that the random-field, long-range Ising model for $\alpha > d+1$, for dimensions $d \geq 3$, undergoes a phase transition, by making use of contours of the form,

\begin{align*}
     \Gamma_0 \big( n \big) \equiv \big\{ \text{paths }  \gamma \in \Gamma : 0 \in I \big( \gamma \big) , \big| \gamma \big| = n \big\}     \text{, }  
\end{align*}

\noindent which denotes each possible contour $\gamma$, of length $n$, in which the interior of the contour contains the origin $0$, and is of length $n$, which are the maximal connected components of the union of faces $C_x \cap C_y$, for which $\sigma_x \neq \sigma_y$ from the set of all possible contours $\Gamma$. Within each $\gamma$, the \textit{Perierls' type argument} entails reversing the direction of the spins contained within the contour, ie flipping the spins to $-\sigma_i$ and otherwise setting all of the spins outside of the contour as $\sigma_i$. The procedure of performing the spin flipping within the interiors of suitable contours allows for the phase transition of RFIM, in lower dimensions, to be shown to hold. Under the action of performing the spin flip, the distribution of $+$, or $-$, spins within the interior of the contour can be formulated in terms of the difference between the original RFIM Hamiltonian with the RFIM Hamiltonian with countably many spins flipped.

Denoting the spin transformation for suitable contours within the system as $\big( \tau_A \big( \sigma \big) \big)_i : \textbf{R}^{\textbf{Z}^d} \longrightarrow \textbf{R}^{\textbf{Z}^d}$, the mapping takes the form,

\[
\big( \tau_A \big( \sigma \big) \big)_i  \equiv  \text{ } 
\left\{\!\begin{array}{ll@{}>{{}}l}    - \sigma_i       &  \text{, if } 
 i \in A \text{, }  \\
\sigma_i  &    \text{ otherwise } \text{ . }  \\
\end{array}\right.
\]

\noindent With $\Gamma$, $\Gamma_0 \big( n \big)$ and $\big( \tau_A \big( \sigma \big) \big)_i \text{ }$, a portion of previous results for demonstrating that the phase transition occurs for the random-field, long-range Ising model are captured with the following \textbf{Proposition}.

\bigskip

\noindent \textbf{Proposition} \textit{1} (\textit{the impact of reversing spins inside contours for the long range, random field Ising model Hamiltonian under plus boundary conditions}, [1], \textbf{Proposition} \textit{2.1}). For $\alpha > d+1$, there exists a constant $c>0$ such that, for the random-field, long-range Ising model at inverse temperature $\beta > 0$,

\begin{align*}
 \mathcal{H}^{\mathrm{LR-RF}, +}_{\Lambda} \big( \tau_{\gamma} \big( \sigma \big)  \big) -  \mathcal{H}^{\mathrm{LR-RF}, +}_{\Lambda} \big(  \sigma \big)  \leq - J c \big| \gamma \big|   \text{. } 
\end{align*}

\noindent The \textbf{Proposition} above demonstrates the impact of reversing the spins contained within $\gamma$, under $\tau_{\gamma} \big( \sigma \big)$, with the spins $\sigma$ before $\tau_{\gamma} \big( \cdot \big)$ is applied. Along similar lines, from the density introduced previously under $+$ boundary conditions with $\mathcal{D}^{\mathrm{LR-RF}+}_{\Lambda , \beta} \big( \cdot , \cdot \big) \equiv \mathcal{D}^{+}_{\Lambda , \beta} \big( \cdot , \cdot \big)$, the equality,

\begin{align*}
      \frac{\mathcal{D}^{+}_{\Lambda , \beta} \big( \sigma , h \big)   Z^{+ }_{\Lambda ,\beta }  \big( h \big) }{\mathcal{D}^{+}_{\Lambda , \beta} \big( \tau_{\gamma} \big( \sigma \big)  , \tau_{\gamma} \big( h \big) \big)   Z^{+}_{\Lambda ,\beta }  \big( \tau \big( h \big) \big)  }    =  \mathrm{exp} \big[      \beta    \mathcal{H}^{\mathrm{LR-RF}, +}_{\Lambda} \big( \tau_{\gamma} \big( \sigma \big)  \big) - \beta  \mathcal{H}^{\mathrm{LR-RF}, +}_{\Lambda} \big(  \sigma      \big)                 \big] \text{ }     \text{, }  
\end{align*}

\noindent between the ratio of the product of the density $\mathcal{D}^{+}_{\Lambda , \beta} \big( \sigma , h \big)$, and $ Z^{+ }_{\Lambda ,\beta }  \big( h \big) $, with the product of $\mathcal{D}^{+}_{\Lambda , \beta} \big( \tau_{\gamma} \big( \sigma \big)  , \tau_{\gamma} \big( h \big) \big)$, and $ Z^{\eta}_{\Lambda ,\beta }  \big( \tau \big( h \big) \big)$, is equivalent to the exponential of the difference between the long-range, random-field Ising model under $\tau_{\gamma} \big( \sigma \big)$ and $\sigma$, respectively. The exponential of the difference of the original, and modified, Hamiltonians equaling the ratio of the product of an exponential of the long-range probability measure with the partition function is of great significance from the Peierls' argument, as the impact of the spin-flipping procedures within interior of countours distinguishes probability measures under $+$ or $-$ boundary conditions. Similarly, under the probability measure and distributions functions for the long range Ising model, instead for $\mathcal{D}^{\mathrm{LR}+}_{\Lambda , \beta} \big( \cdot , \cdot \big) \equiv \mathcal{D}^{+}_{\Lambda , \beta} \big( \cdot , \cdot \big)$,

\begin{align*}
\frac{\mathcal{D}^{+}_{\Lambda , \beta} \big( \sigma , \eta \big)   Z^{+ }_{\Lambda ,\beta }  \big( \eta \big) }{\mathcal{D}^{+}_{\Lambda , \beta} \big( \tau_{\gamma} \big( \sigma \big)  , \tau_{\gamma} \big( \eta \big) \big)   Z^{+}_{\Lambda ,\beta }  \big( \tau \big( \eta \big) \big)  }    =  \mathrm{exp} \big[      \beta    \mathcal{H}^{\mathrm{LR}, +}_{\Lambda} \big( \tau_{\gamma} \big( \sigma \big)  \big) - \beta  \mathcal{H}^{\mathrm{LR}, +}_{\Lambda} \big(  \sigma      \big)                 \big] \text{ }     \text{. } 
\end{align*}

\noindent The equality between the exponential of the difference of two long-range Hamiltonians

Under a random, external field introduced with iid, Gaussian $\big\{ h_x \big\}$, it is possible for the partition function $Z^{+}_{\Lambda , \beta} \big( \tau \big( h \big) \big)$ to exceed $Z^{+}_{\Lambda , \beta} \big( h \big)$. If this were the case, the parameter,

\begin{align*}
   \Delta_A \big( h \big) \equiv - \frac{1}{\beta} \mathrm{log} \bigg[   \frac{Z^{+}_{\Lambda, \beta} \big( h \big) }{Z^{+}_{\Lambda, \beta} \big( \tau_A \big( h \big) \big) }     \bigg] \text{, }  
\end{align*}

\noindent captures the probability of such an event occurring, in which there exists a path, sampled from $\Gamma$, for which,

\begin{align*}
  \underset{\gamma \in \Gamma_0 }{\mathrm{sup}} \frac{\big|  \Delta_{I ( \gamma )} \big( h \big) \big|}{c_1 \big|  \gamma \big| }  < \frac{1}{4} \text{, }  
\end{align*}

\noindent which we denote with the 'bad' event, $\mathcal{B}$. Hence the complementary event for a bad event is given by,

\begin{align*}
 \mathcal{B}^c \equiv \bigg\{  \underset{\gamma \in \Gamma_0 }{\mathrm{sup}} \frac{\big|  \Delta_{I ( \gamma )} \big( h \big) \big|}{c_1 \big|  \gamma \big| }  > \frac{1}{4} \bigg\}    \text{. } 
\end{align*}

\noindent From the supremum introduced above, of a term inversely proportional to the length, and directly proportional to the interior of each such $\gamma$, several bounds leading up to the \textit{Peierls' argument} incorporate $\tau_A \big( \sigma \big)$, one of which is first introduced below. From the probability measures $\textbf{P}^{\mathrm{LR-RF}}_{\Lambda} \big( \cdot \big)$, and $\textbf{P}^{\mathrm{LR}}_{\Lambda} \big( \cdot \big)$, denote $\textbf{P}^{\mathrm{RF}}_{\Lambda} \big( \cdot \big)$ as the probability measure for the random field Ising model.

\bigskip

\noindent \textbf{Lemma} \textit{1} (\textit{constant times an exponential upper bound for the random field Ising model} [1], \textbf{Lemma} \textit{3.4}). For $A , A^{\prime} \subsetneq \textbf{Z}^d$, with $A \cap A^{\prime} \neq \emptyset$ and $\big| A \big| , \big| A^{\prime} \big| < + \infty$,

\begin{align*}
  \textbf{P}^{\mathrm{RF},+}_{\Lambda} \big[   \big|  \Delta_A \big( h \big)       \big| \geq \lambda \big| h_{A^c}      \big]  \leq 2 \text{ } \mathrm{exp} \bigg[        - \frac{\lambda^2}{8 e^2 \big| A \big|}           \bigg]  \text{, }  
\end{align*}

\noindent and also that,

\begin{align*}
    \textbf{P}^{\mathrm{RF},+}_{\Lambda} \big[   \big|  \Delta_A \big( h \big) - \Delta_{A^{\prime}} \big( h \big)    \big| > \lambda \big|   h_{(A \cup A^{\prime})^c}     \big]   \leq 2 \text{ } \mathrm{exp} \bigg[     - \frac{\lambda^2}{8 e^2 \big| A \Delta A^{\prime} \big|}                \bigg]   \text{, }  
\end{align*}

\noindent for the symmetric difference between the sets $A$ and $A^{\prime}$, $A \Delta A^{\prime}$.

\bigskip

\noindent Besides the result above, we must also make use of a coarse-graining procedure. For the procedure, as described in [1] and [2], introduce a coarse grained renormalization of $\textbf{Z}^d$,

\begin{align*}
 C_m \big( x \big) \equiv      \bigg[        \overset{d}{\underset{i=1}{\prod}}  \big[ 2^m x_i - 2^{m-1} , 2^m x_i + 2^{m-1}     \big]  \bigg] \cap \textbf{Z}^d        \text{, }  
\end{align*}

\noindent corresponding to the cube over the hypercube, with center at $2^m x$, with side length $2^m-1$, an \textit{m-cube}, which is a restatement of the coarse-graining approach of [5]. From the object above, we make use of the convention that $C_0 \big( 0 \big)$ denotes the point about $0$. Additionally, denote,

\begin{align*}
       \mathcal{P}_i        \big(   A \cap \mathcal{R} \big) \equiv \big\{ x \in \mathcal{R}_i          :    l^i_x \cap A \neq \emptyset \big\}     \text{, }  
\end{align*}

\noindent which also satisfies,

    \begin{align*}
 \mathcal{P}_i \big( A \cap \mathcal{R} \big) \supsetneq \underset{1 \leq i \leq d}{\bigcup} \big(  \mathcal{P}^{\mathrm{G}}_i \big( A \cap \mathcal{R} \big)  \cup  
\mathcal{P}^{\mathrm{B}}_i \big( A \cap \mathcal{R} \big) \big)   \text{, }  
\end{align*}

\noindent for a rectangle $\mathcal{R} \equiv \overset{n}{\underset{i=1}{\prod}} \big[ 1 , r_i \big]$, with $\mathcal{R}_i \cap \big[ 1 , r_i \big] \neq \emptyset$ for every $i$, which is given by,

\begin{align*}
\mathcal{R} \supsetneq \mathcal{R}_i \equiv \big\{ x \in \mathcal{R} : x_i = 1 \big\}  \text{ } \text{ ,} 
\end{align*}

\noindent and $l^i_x \equiv \big\{       x + k e_i   :   1 \leq k \leq r_i \big\}$, satisfying $\mathcal{R} \cap \mathcal{R}_i \neq \emptyset$ for each $i$, as the set of points for which $l^i_x \cap A \neq \emptyset$. From this, denote the \textit{good} set of points in the plane,

\begin{align*}
     \mathcal{P}^{\mathrm{G}}_i \big( A \cap \mathcal{R} \big)  \equiv         \big\{   \forall \text{ rectangles } \mathcal{R}_i \text{ , } \exists \text{ } \text{countably many }     x \in \mathcal{P}_i \big( A \cap \mathcal{R}  \big) : l^i_x \cap \big( A \backslash \mathcal{R}  \big) \neq \emptyset     \big\}   \text{, }  
\end{align*}

\noindent and, similarly, denote the set of bad points,

\begin{align*}
          \mathcal{P}^{\mathrm{B}}_i \big( A \cap \mathcal{R} \big)  \equiv  \big( \mathcal{P}^{\mathrm{G}} \big( A \cap \mathcal{R} \big) \big)^c              \text{, }  
\end{align*}

\noindent for which $l^i_x \cap \big( A \backslash \mathcal{R}  \big)  \equiv \emptyset$. In comparison to the contours discussed in [2] which are used to implement a \textit{Peirels' argument}, related to the projections $\mathcal{P}_i$, that,

\begin{align*}
    \big| \mathcal{P}^{\mathrm{G}}_i \big( A \cap \mathcal{R} \big) \big|      \leq       \big|    \partial_{\mathrm{ex}} A \cap \mathcal{R}                \big|  \text{, }  
\end{align*}

\noindent in which the set of \textit{good} points has cardinality less than, or equal to, the cardinality of $\partial_{\mathrm{ex}} A \cap \mathcal{R}$, where,

\begin{align*}
  \partial_{\mathrm{ex}} A  \equiv  \big\{ \forall v \in A^c \cup \partial A , \exists \text{ }  v^{\prime} \in \partial A : v \cap v^{\prime} \neq \emptyset  \big\}  \text{, }  
\end{align*}

\noindent and,

\begin{align*}
        \big|  \mathcal{P}^{\mathrm{B}}_i \big( A \cap \mathcal{R} \big)     \big|      \leq C       \big|    \mathcal{R}_d            \big|    \text{, }  
\end{align*}

\noindent in which the set of \textit{bad} points has cardinality less than, or equal to, the cardinality of a rectangular subset of the hypercube, $\mathcal{R}_d$, for a real parameter $C \equiv \frac{\lambda}{r_j}$, while finally, that,

\begin{align*}
      \overset{d}{ \underset{i=1}{\sum}}   \big| \mathcal{P}_i \big( A \cap \mathcal{R} \big) \big|      \leq    c    \big|    \partial_{\mathrm{ex}} A \cap \mathcal{R}      \big|      \text{, }  
\end{align*}

\noindent where the exterior boundary of a path is given by,

\begin{align*}
      \partial_{\mathrm{ex}} \big( \Lambda \big) \equiv \big\{  \forall   x \in \Lambda^c    \text{ } , \text{ } \exists y \in \Lambda : \big| x - y \big| = 1   \big\}          \text{. } 
\end{align*}

\noindent Similarly, the interior boundary of a path is given by,

\begin{align*}
          \partial_{\mathrm{int}}  \big( \Lambda \big) \equiv \big\{   \forall x \in \Lambda \text{ } , \text{ } \exists y \in \Lambda^c : \big| x - y \big| = 1   \big\}              \text{. } 
\end{align*}

\noindent Above, the summation of the cardinality of the set of \textit{all} points in the projection $\mathcal{P}_i$ is less than, or equal to, $\partial_{\mathrm{ex}} A \cap \mathcal{R}$, for every $1 \leq i \leq d$, and some $c>0$. Following a description of the paper organization in the next section, we distinguish between the types of contours discussed in [1], and in [2].

\subsection{Paper organization}

\noindent With the definition of the long range, random-field, and long range, random-field Ising models, in the next section we differentiate between contours discussed in [1] and [2], from which the existence of a phase transition can be provided for the long range, random-field Ising model for $d < \alpha \leq d+1$. In order to adapt the argument provided in [1] with the contours described in [2], we implement several steps of the argument for the long range contour system surrounding the coarse graining procedure. Fundamentally, the coarse graining procedure introduces a partition of the square lattice; through such a partition, the discontinuity of the phase transition is a statement surrounding the existence of more than two extremal probability measures, or otherwise referred to as Gibbs states. Altogether, to exhibit that a phase transition occurs for lower dimensions in the long range Ising model, we prove the following result.

\bigskip

\noindent \textbf{Theorem PT} (\textit{the long range, random-field Ising model undergoes a phase transition in lower dimensions}). Over a finite volume $\Lambda$, for $d \geq 3$, there exists a critical parameter $\beta_c$, with $\beta_c \equiv \beta_c \big( \alpha , d \big)$, and another parameter $\epsilon$, with $\epsilon \equiv \epsilon \big( \alpha , d \big)$, so that for parameters $\beta \geq \beta_c$ and $\epsilon \leq \epsilon_c$,

\begin{align*}
  \textbf{P}^{\mathrm{LR},+}_{\Lambda,\beta, \epsilon} \neq \textbf{P}^{\mathrm{LR},-}_{\Lambda,\beta, \epsilon} \text{, }  
\end{align*}

\noindent $\textbf{P}$-almost surely, in which the long range measures under $+$ and $-$ boundary conditions are not equal.

\section{Contours in the long range, random-field Ising model for the Peirels' argument}

\noindent We introduce long range contours below. As individual paths that are part of a larger collection of suitable contours, the interior of spins within the interior being reversed distinguishes the long-range probability measures under $+$, and $-$, boundary conditions. As such, the phase transitions for lower, and higher, dimensional RFIMs depends upon the number of suitable contours that one can perform the spin flipping procedure on, in addition to an upper bound for $ \mathcal{H}^{\mathrm{LR}, +}_{\Lambda} \big(  \big( \tau_{\Gamma} \big( \sigma \big)  \big)_x \big)  -  \mathcal{H}^{\mathrm{LR}, +}_{\Lambda} \big(  \sigma \big)     $. The upper bound for the later item, besides previously appearing in the expression,

\begin{align*}
      \frac{\mathcal{D}^{+}_{\Lambda , \beta} \big( \sigma , h \big)   Z^{+ }_{\Lambda ,\beta }  \big( h \big) }{\mathcal{D}^{+}_{\Lambda , \beta} \big( \tau_{\gamma} \big( \sigma \big)  , \tau_{\gamma} \big( h \big) \big)   Z^{+}_{\Lambda ,\beta }  \big( \tau \big( h \big) \big)  }    =  \mathrm{exp} \big[      \beta    \mathcal{H}^{\mathrm{LR-RF}, +}_{\Lambda} \big( \tau_{\gamma} \big( \sigma \big)  \big) - \beta  \mathcal{H}^{\mathrm{LR-RF}, +}_{\Lambda} \big(  \sigma      \big)                 \big] \text{ }     \text{, }  
\end{align*}

\noindent is also of significance for upper bounding the probability,

\begin{align*}
   \textbf{P} \bigg[   \underset{\gamma \in \Gamma_0}{\mathrm{sup}} \frac{\big|  \Delta_{I_{-} ( \gamma )} \big( \eta \big) \big|}{ \big|  \gamma \big| }  > 1     \bigg] \leq \mathrm{exp} \big( - C^{\prime}_2 \epsilon^{-2} \big)     \text{. } 
\end{align*}

\noindent The probability above determines whether the length of suitable contours, as a normalization to the total number of contours belonging to $\Delta_{I_{-} ( \gamma )} \big( \eta \big)$, exceeds $1$. In demonstrating that the above probability is exponentially small, Dudley's argument is used to upper bound the expectation,

\begin{align*}
   \textbf{E} \bigg[   \text{ } \underset{_{\gamma \in \Gamma_0 ( n ) }   }{\mathrm{sup}}\Delta_{I_{-} ( \gamma )} \big( \eta \big) \text{ }     \bigg] \text{. } \end{align*}

\subsection{Contours for the long range Ising model}

\noindent To introduce another family of contours for the \textit{Peierls’ argument}, consider the following.

\bigskip

\noindent \textbf{Definition} \textit{1} (\textit{new contours for the Peierls’ argument}, [2]). For the long range Ising model, real $M,a,r >0$, and a configuration $\sigma \in \Omega^{\mathrm{LR}}$, the sample space of all long range Ising model configurations, from the boundary $\partial \sigma$, the set of all $\big( M , a , r \big)$-partitions, $\Gamma \big( \sigma \big) \equiv \big\{ \bar{\gamma} : \bar{\gamma} \subset \partial \sigma \big\} \neq \emptyset$, satisfies:

\begin{itemize}
    \item[$\bullet$] \textit{Property 1} (\textit{partition equality}): Given $\Gamma \big( \sigma \big)$, there exists countably many $\bar{\gamma}$ which partition each $\partial \sigma$, in which $\underset{\bar{\gamma} \in \Gamma ( \sigma )}{\cup} \bar{\gamma} \equiv \partial \sigma$, such that for another path $\bar{\gamma}^{\prime}$, with $\bar{\gamma} \cap \bar{\gamma}^{\prime} \neq \emptyset$, $\bar{\gamma}^{\prime}$ is contained in the connected component of $( \bar{\gamma} )^c$.

    \item[$\bullet$] \textit{Property 2} (\textit{decomposing each} $\bar{\gamma}$). For all $\bar{\gamma} \in \Gamma \big( \sigma \big)$, $\exists $ $1 \leq n \leq 2^r - 1$ such that:

        \begin{itemize}
            \item[$\bullet$] \textit{Property 2A}: $\bar{\gamma}$ can be expressed with the union $\bar{\gamma} \equiv \underset{1 \leq k \leq n}{\bigcup} \bar{\gamma}_k$, for $\bar{\gamma}_k$ such that $\bar{\gamma}_k \cap \bar{\gamma} \neq \emptyset$ for every $k$.

            \item[$\bullet$] \textit{Property 2B}: For $\bar{\gamma}, \bar{\gamma}^{\prime} \in \Gamma \big( \sigma \big)$ such that $\bar{\gamma} \cap \bar{\gamma}^{\prime} \neq \emptyset$, there exists two strictly positive $n \neq n^{\prime}$, for which,

            \begin{align*}
              \mathrm{d} \big(  \bar{\gamma} , \bar{\gamma}^{\prime} \big) > M \text{ } \mathrm{min} \big\{   \underset{1 \leq k \leq n}{\mathrm{max}}    \mathrm{diam} \big( \bar{\gamma}_k \big)          , \underset{1 \leq j \leq n^{\prime}}{\mathrm{max}}\mathrm{diam} \big( \bar{\gamma}^{\prime}_j \big)                 \big\}^a  \text{, } 
            \end{align*}

            with respect to the metric $\mathrm{d} \big( \cdot , \cdot \big)$ between paths belonging to $\Gamma \big( \sigma \big)$, where,

            \begin{align*}
                \mathrm{d} \big( \gamma_1 , \gamma_2 \big) \equiv \big\{ \forall  n \in \textbf{Z}_{\geq 0} \text{, }  \text{ } \exists \text{ } \gamma_1 , \gamma_2 \in \Gamma :   \big\| \gamma_1 - \gamma_2 \big\|_1 = n     \big\}              \text{. } 
            \end{align*}

            \end{itemize}
    
\end{itemize}

\noindent With \textbf{Definition} \textit{1}, we also denote the set of all \textit{connected components} of some $\sigma$ in finite volume, below.

\bigskip

\noindent \textbf{Definition} \textit{2} (\textit{connected components in a finite volume}). For any $m_1 \neq m_2 >0$, and two vertices $x \neq x^{\prime}$, there exists two \textit{m-cubes}, $C_{m_1} \big( x \big)$ and $C_{m_2} \big( x^{\prime} \big)$, such that the edge set,

\begin{align*}
V_n \equiv v \big( G_n \big( \Lambda \big) \big)  \equiv \big\{  v \in C_m \big( x \big) : v \cap V \big( \Lambda \big)  \neq   \emptyset  \big\}    \text{, }  
\end{align*}

\noindent is comprised of the minimum number of cubes for which the union of \textit{m-cubes} covers the set of \textit{connected components}, while the \textit{edge set},

\begin{align*}
 E_n \equiv e \big( G_n \big( \Lambda \big) \big) \equiv  \big\{    e \in E \big( \Lambda \big) :  \big| e \cap  E \big( \Lambda \big)  \cap C_m \big( x \big) \big|  \leq M d^a 2^{a n} \big\}    \text{, }  
\end{align*}

\noindent is comprised of the number of edges that have nonempty intersection with $E \big( \Lambda \big)$ and $C_m \big( x \big)$, for $G_n \big( \Lambda \big) \equiv \big(  V_G , E_G \big)$. Denote the set of \textit{connected components}, $\mathscr{G}_n \big( \Lambda\big)$, associated with some configuration, and contained with some \textit{m-cube}, as,

\begin{align*}
    \gamma_{G} \big( \Lambda , C_m \big( x \big)  \big) \equiv \gamma_G \equiv \underset{G_i \cap \Lambda \cap C_m ( x) \neq \emptyset}{\underset{G_i \subsetneq G}{\bigcup}} \gamma_{G_i} \equiv \underset{\forall C_m ( x ) v \in V_G : C_m ( x ) \cap   v \neq \emptyset   }{\bigcup}   \big(  \Lambda  \cap  C_m \big( x \big)      \big)     \text{, }  
\end{align*}

\noindent corresponding to the \textit{connected components} with nonempty intersection with an \textit{m-cube}.

\bigskip

\noindent With the set of \textit{connected components} from \textbf{Definition} \textit{2}, denote a set of partitions, $\big\{ \mathscr{P}_i \big\}_{i \in \mathcal{I}}$ for some countable index set $\mathcal{I}$, such that $\mathscr{P}_i \cap G_n \big( \Lambda \big) \neq \emptyset$ for every $i$, as the set of finite subvolumes of $\Lambda$ for which,

\[
\mathscr{P}_i \equiv  \text{ } 
\left\{\!\begin{array}{ll@{}>{{}}l}    \big\{ \forall   G \in   \mathscr{G}_n \big( \Lambda \big) , \exists \sigma_i , r > 0 :  \mathscr{G}_n \big( \sigma_i  \big) \cap \Lambda \neq \emptyset ,  \big| v \big( G \big)  \big|  \leq 2^r - 1   \big\}      &  \text{, if } i > 0
 \text{, }  \\ \big\{ \forall   G \in   \mathscr{G}_n \big( \Lambda \big) , \exists \sigma_i , r > 0  : \mathscr{G}_n \big( \sigma_i  \big) \cap \Lambda \neq \emptyset ,  1 \leq \big| v \big( G \big)  \big|  \leq 2^r - 1   \big\}     
  &    \text{, if } i \equiv 0 \text{ }  \text{ . }  \\
\end{array}\right.
\]

\noindent $\mathscr{P}_i$ is otherwise assumed to be equal to $\emptyset$ if $\partial \sigma_i = \emptyset$. From \textbf{Proposition} \textit{3.5} in [2], the collection $\big\{ \mathscr{P}_i \big\}$ satisfies \underline{Property 1}, and \underline{Property 2}. Finally, below, introduce the \textit{inner boundary} and the set of edges that are exactly incident with the boundary configuration.

\bigskip

\noindent \textbf{Definition} \textit{3} (\textit{inner and incident boundaries of edges to the boundary configuration}). Denote the \textit{inner boundary of edges} to $\partial \sigma_i$ with,

\begin{align*}
    \partial_{\mathrm{in}}  \big( \Lambda  ,  \partial \sigma_i \big)  \equiv \partial_{\mathrm{in}} \Lambda   \equiv  \big\{  \forall  \sigma_i , \exists m > 0  :  \big(  \mathscr{G}_n \big( \Lambda \big)    \cap  C_m \big( x \big)   \big) \cap \partial \sigma_i  \equiv \emptyset       \big\}  \text{, }  
\end{align*}

\noindent and the \textit{incident boundary of edges to} $\partial \sigma_i$ with,

\begin{align*}
     \mathcal{B} \big( \partial \sigma_i \big)    \equiv \big\{ \forall  \sigma_i , \exists m > 0 :     \big| \mathscr{G}_n \big( \Lambda \big) \cap    C_m \big( x \big)  \big| \equiv \big| \partial \big( \mathscr{G}_n \big( \Lambda \big) \big) \big|         \big\}  \text{, }  
\end{align*}

\noindent under the assumption that $\partial_{\mathrm{in}} \Lambda, \mathcal{B} \big( \partial \sigma_i \big) \neq \emptyset$.

\bigskip

\noindent From quantities from \textbf{Definiton} \textit{3}, the isoperimetric inequality states,

\begin{align*}
     \big|   \Lambda   \big|^{1 - \frac{1}{d}}        \leq   \big| \partial_{\mathrm{in}} \Lambda \big|  \text{, }  
\end{align*}

\noindent for the dimension $d$.

\subsection{Long range, versus long range, random-field Ising model contours}

\noindent From contours for the long range Ising model of the previous section, the procedure for reversing the orientation of spins differs. Below, we introduce essential components of the coarse graining procedure, which not only determines the average number of spins that are flipped within the interior of suitable contours, but also controls the discrepancy between the long-range probability measure under $+$, and $-$, boundary conditions. The fact that the statement of the phase transition for the lower dimensional RFIM is dependent upon probability measures under $+$, and $-$, boundary conditions implies that the transformation of reversing spins within the interior of suitable contours can significantly impact the probability of an event occurring. To this end, another important component of the argument relies upon the statement that the probability, supported over some $\Delta$, of the complement of bad events, $\mathcal{B}$, occurring is exponentially small,

\begin{align*}
    \textbf{P}_{\Lambda} \big[ \mathcal{B}^c \big]    \leq \mathrm{exp} \big( - C_1 \epsilon^{-2} \big)        \text{, } 
\end{align*}

\noindent given $C_1 \neq \epsilon$ taken to be sufficiently small. Fix the \textit{m-cube} of side length $m$ about the point $0$,

\begin{align*}
  C_0 \big( m \big) \equiv \big\{ \mathrm{sp} \big( \gamma \big) \subsetneq \textbf{Z}^d , \big| \mathrm{sp} \big( \gamma \big) \big| < + \infty : \gamma \in \mathcal{E}^{-}_{\Lambda} , 0 \in V \big( \gamma \big) , \big| \gamma \big| = m \big\}   \text{. } 
\end{align*}

\noindent As opposed to $\big( \tau_A \big( \sigma \big) \big)_i$ for countours in the long range, random-field Ising model, the flipping procedure is, for the set $\Gamma$ at each $x$, given by the map $\big( \tau_{\Gamma} \big( \sigma \big) \big)_x:  \Omega \big( \Gamma \big)    \longrightarrow \Omega^{-}_{\Lambda}$, where the target space of the mapping is,

\begin{align*}
  \Omega^{-}_{\Lambda} = \big\{           \text{collection of all paths contained in } \Lambda \text{ with -} 1 \text{ labels}    \big\}  \equiv  \big\{ \gamma \in \Lambda : \gamma \cap \Lambda \neq \emptyset , \\   \mathrm{lab} \big( \gamma \big) \equiv - 1  \big\}  \text{, }  
\end{align*}

\noindent as,

\[
\big( \tau^{\mathrm{LR}}_{\Gamma} \big( \sigma \big) \big)_x \equiv \big(\tau^{\mathrm{LR}} \big( \sigma \big)\big)_x \equiv \big( \tau_{\Gamma} \big( \sigma \big) \big)_x  \equiv  \text{ } 
\left\{\!\begin{array}{ll@{}>{{}}l}     \sigma_x   &  \text{, if } 
 x \in I_{-} \big( \Gamma \big) \cup V \big( \Gamma \big)^c \text{, }  \\
- \sigma_x  &    \text{, if } x \in I_{+} \big( \Gamma \big)  \text{                    } \text{ , }  \\   - 1  &   \text{      
  , if } x \in \mathrm{sp} \big(\Gamma \big) \text{, }  \\ 
\end{array}\right.
\]

\noindent which can be expressed with the following over all $n$ components of $\gamma$, with, 

\begin{align*}
      \big( \tau_{\Gamma} \big( \sigma \big) \big)_x =  \big( \tau_{ \{ \gamma_1 , \cdots , \gamma_n \} } \big( \sigma \big) \big)_x   \text{. } 
\end{align*}

\noindent Also, given the support, collection of edges with $-$ labels, the set of all labels, vertices of $G$, and interior of each $\gamma$, each of which are respectively given by,

\begin{align*}
   \big| \gamma \big| \equiv  \mathrm{sp} \big( \gamma \big) \equiv \big\{   \text{support of paths }  \gamma     \big\}   \text{, }  \\  \mathcal{E}^{-}_{\Lambda} \equiv \big\{   \forall     \Gamma \equiv \big\{ \gamma_1 , \cdots , \gamma_n \big\} \text{ } ,  \text{ } \exists V \big( \Gamma \big) \subset \Lambda  : \textit{compatible} \text{ }  \Gamma , \textit{external}  \text{ } \gamma_i , \mathrm{lab} \big( \gamma_i \big) = - 1       \big\}  \text{, }   \\ \mathrm{lab}_{\bar{\gamma}} \equiv \big\{ \text{labels of paths } \gamma \big\} \equiv \underset{n\geq 0}{\underset{\text{paths } \gamma}{\bigcup} } \big\{ \forall  i >  0  , \bar{\gamma} \equiv \big( \bar{\gamma}^0 , \cdots , \bar{\gamma}^n \big)  \in \Gamma , \exists 1 < i < n   :  \bar{\gamma}^i \longrightarrow \big\{ - 1 , + 1 \big\} \big\}   \text{, }  \\  V \big( G \big) \supsetneq V \big( \Gamma \big) \equiv \big\{   v \in  v \big( G \big) :  v \cap G \cap \Lambda  \neq \emptyset      \big\}  \text{, }  
 \\ I_{\pm} \big( \gamma \big)  \equiv \underset{k \geq 1, 1 \leq k \leq n}{\bigcup} I_{\pm} \big( \gamma_k \big)  \equiv \underset{\mathrm{lab}_{\bar{\gamma}} (  I  ) = \pm 1 }{\underset{k\geq 1}{\bigcup}} I \big( \mathrm{sp} \big( \gamma \big) \big)^{k}  \text{, }  
\end{align*}

\noindent in addition to the two quantities,

\begin{align*}
  V \big( \gamma \big) \equiv \mathrm{sp} \big( \gamma \big) \cup I \big( \gamma \big) \equiv \mathrm{sp} \big( \gamma \big) \cup      \underbrace{  \big( 
 {I_{+} \big( \gamma \big) \cup I_{-} \big( \gamma \big)}   \big) }_{I ( \gamma ) \equiv I_{+} ( \gamma ) \cup I_{-} ( \gamma )}   \text{, }  
\end{align*}

\noindent where in the definition of $\mathcal{E}^{-}_{\Lambda}$, paths are considered \textit{compatible} from the set of all paths $\Gamma$ if there exists a configuration from the long range sample space, $\sigma$, whose contours coincide with those of $\Gamma$. Similarly, for paths with $+1$ labels, introduce the collection of \textit{compatible} paths over $\Lambda$,

\begin{align*}
  \mathcal{E}^{+}_{\Lambda} \equiv \big\{   \forall     \Gamma \equiv \big\{ \gamma_1 , \cdots , \gamma_n \big\}, \text{ } \exists V \big( \Gamma \big) \subset \Lambda  : \textit{compatible} \text{ }  \Gamma , \textit{external}  \text{ } \gamma_i , \mathrm{lab} \big( \gamma_i \big) = + 1       \big\}  \text{, }   \text{. } 
\end{align*}

\noindent From the quantities introduced above that are assocated with the flipping procedure $\big( \tau_{\Gamma} \big( \sigma \big) \big)_x$, it is also important to state the difference in $\mathcal{H}^{\mathrm{LR-RF}, +}_{\Lambda} \big( \tau_{\gamma} \big( \sigma \big)  \big) -  \mathcal{H}^{\mathrm{LR-RF}, +}_{\Lambda} \big(  \sigma \big)$ between $\tau_{\gamma} \big( \sigma \big)$ and $\sigma$. For the long range Ising model with the contour system defined in \textit{2.1}, the long range Hamiltonian instead satisfies, (\textbf{Proposition} \textit{4.5}, [2]),

\begin{align*}
\mathcal{H}^{\mathrm{LR},-}_{\Lambda} \big( \tau \big( \sigma \big)  \big) -    \mathcal{H}^{\mathrm{LR},-}_{\Lambda} \big( \sigma \big)    \leq   -  c_1 \big| \gamma \big| - c_2 F_{I_{+}( \gamma) } - c_3 F_{\mathrm{sp}(\gamma)}  \text{, }  
\end{align*}

\noindent for a long range configuration $\sigma$, strictly positive $c_1, c_2, c_3$, and for the functions,

\begin{align*}
        F_{I_{\pm} ( \gamma) }  \equiv  \underset{y \in (I_{\pm}( \gamma))^c}{\underset{x \in I_{\pm} ( \gamma)}{\sum}} J_{x,y} \text{, }  \\ F_{\mathrm{sp}(\gamma)} \equiv \underset{y \in (\mathrm{sp}(\gamma))^c}{\underset{x \in \mathrm{sp}(\gamma)}{\sum}} J_{x,y}  \text{. } 
\end{align*}  

\bigskip

\noindent Long range contours differ from long range, random-field contours to a similar condition as raised in the isoperimetric inequality, in which, (\textbf{Lemma} \textit{4.3}, [2]),

\begin{align*}
  \mathrm{diam} \big( \Lambda \big) \geq k_d \big| \Lambda \big|^{\frac{1}{d}}         \text{, }  
\end{align*}

\noindent in which the diameter of each such path is bound below by some strictly positive prefactor times the cardinality of the finite volume, $\Lambda$, in addition to the fact that the paths for the long range Ising model, in comparison to those from the long range, random-field Ising model, do not satisfy,

\begin{align*}
    \mathscr{C}_l \big( \gamma \big) \equiv \underset{l \in \textbf{N}}{\bigcup}    \big\{ C_l : \big| C_l \cap I \big( \gamma \big) \big| \geq \frac{1}{2} \big| C_l \big| \big\}    \text{, }  
\end{align*}

\noindent introduced as the $C_l$ admissibility condition [1], which has boundary,

\begin{align*}
  \partial  \mathscr{C}_l \big( \gamma \big) \equiv \big\{     \big( C_l , C^{\prime}_l \big) : C^{\prime}_l \not\in \mathscr{C}_l \big( \gamma \big) , \big| C^{\prime}_l \cap C_l \big| = 1 \big\}   \text{. } 
\end{align*}

\section{Phase transition for the long-range Ising model}

\noindent The argument for proving that a phase transition occurs for the long range, random field Ising model can be applied to demonstrate that a phase transition occurs for the long range Ising model, beginning with the following.

\subsection{Beginning the argument}

\noindent We must determine the upper bound for the behavior of the long range Ising model Hamiltonian under the flipping procedure given in the previous section with $\big( \tau^{\mathrm{LR}} \big( \sigma \big) \big)_x$. For a new range of parameters $\alpha$ satisfying $d < \alpha \leq d+1$, instead of upper bounding the difference $\mathcal{H}^{\mathrm{LR}, }_{\Lambda} \big(  \big( \tau_{\Gamma} \big( \sigma \big)  \big)_x \big)  -  \mathcal{H}^{\mathrm{LR}, -}_{\Lambda} \big(  \sigma \big)$, under $-$ boundary conditions in the $\alpha > d+1$ regime, we upper bound the difference $\mathcal{H}^{\mathrm{LR}, +}_{\Lambda} \big(  \big( \tau_{\Gamma} \big( \sigma \big)  \big)_x \big)  -  \mathcal{H}^{\mathrm{LR}, +}_{\Lambda} \big(  \sigma \big)$, under $+$ boundary conditions in the $d < \alpha \leq d+1$ regime.

\bigskip

\noindent \textbf{Proposition} \textit{1} (\textit{upper bound of the flipping procedure of the long range Ising model Hamiltonian with} + \textit{boundary conditions}). For a long range Ising configuration $\sigma \sim \textbf{P}_{\Lambda,\beta}  \big( \cdot , \cdot  \big)$, with energy $\mathcal{H}^{\mathrm{LR}, \eta}_{\Lambda} \big( \sigma)$, the difference of the Hamiltonian under $\big( \tau^{\mathrm{LR}} \big( \sigma \big) \big)_x$ with the Hamiltonian under $\sigma$ satisfies,

\begin{align*}
       \mathcal{H}^{\mathrm{LR}, +}_{\Lambda} \big(  \big( \tau_{\Gamma} \big( \sigma \big)  \big)_x \big)  -  \mathcal{H}^{\mathrm{LR}, +}_{\Lambda} \big(  \sigma \big)       \leq   -  c^{\prime}_1 \big| \gamma \big| - c^{\prime}_2 F_{I_{-}( \gamma) } - c^{\prime}_3 F_{\mathrm{sp}(\gamma)}  \text{, }  
\end{align*}

\bigskip

\noindent for strictly positive $c^{\prime}_1, c^{\prime}_2, c^{\prime}_3$.

\bigskip

\noindent \textit{Proof sketch of Proposition 1}. The argument strongly resembles the strategy used in \textbf{Proposition} \textit{4.5}, [2], in which the authors express each term in the Hamiltonian of the configuration $\sigma$ acted on by the flipping procedure $\big(\tau_{\Gamma} \big( \sigma \big) \big)_x$ for long range contours. Write out the first long range Hamiltonian on the LHS, denoting $\tau_{\Gamma} \big( \sigma_x \big) \equiv \big( \tau_{\Gamma} \big( \sigma \big)_x \big)$, $\tau_{\Gamma} \big( \sigma_y \big) \equiv \big( \tau_{\Gamma} \big( \sigma \big)_y \big)$, $\gamma^1 \equiv \gamma \equiv \{ \gamma^1 , \cdots , \gamma^n \}$, and $\Gamma \big( \sigma \big) \equiv \Gamma$, in which contributions from the Hamiltonian arise from the nonempty regions $I_{-} \big( \gamma \big) \cup V \big( \Gamma \big)^c$, $I_{+} \big( \gamma \big)$, and $\mathrm{sp} \big( \Gamma \big)$, as,

\begin{align*}
    \mathcal{H}^{\mathrm{LR}, +}_{\Lambda} \big( \big( \tau_{\Gamma} \big( \sigma \big)  \big)_x \big) =   - \underset{\Gamma \equiv \cup_i \{ \gamma^i_1 , \cdots , \gamma^i_n \} }{\underset{x,y \in  I_{-} ( \gamma ) \cup V ( \Gamma )^c }{\sum}}   J_{x,y}  \big[ \tau_{\Gamma} \big( \sigma_x \big) \tau_{\Gamma} \big( \sigma_y \big) \big]     - \underset{\Gamma \equiv \cup_i \{ \gamma^i_1 , \cdots , \gamma^i_n \} }{\underset{x,y \in ( I_{-} ( \gamma ) \cup V ( \Gamma )^c )^c }{\sum}}  J_{x,y}  \big[ \tau_{\Gamma} \big( \sigma_x \big) \tau_{\Gamma} \big( \sigma_y \big) \big]       \\ - 
    \underset{\Gamma \equiv \cup_i \{ \gamma^i_1 , \cdots , \gamma^i_n \}}{\underset{ x \in I_{-} ( \gamma ) \cup V ( \Gamma )^c }{\underset{ y \in ( I_{-} ( \gamma ) \cup V ( \Gamma )^c )^c }{\sum}}}  J_{x,y}  \big[ \tau_{\Gamma} \big( \sigma_x \big) \tau_{\Gamma} \big( \sigma_y \big) \big] -  \underset{\Gamma \equiv \cup_i \{ \gamma^i_1 , \cdots , \gamma^i_n \} }{\underset{x,y \in I_{+} ( \gamma )}{\sum}} J_{x,y}  \big[ \tau_{\Gamma} \big( \sigma_x \big) \tau_{\Gamma} \big( \sigma_y \big) \big]  \\ -  \underset{\Gamma \equiv \cup_i \{ \gamma^i_1 , \cdots , \gamma^i_n \} }{\underset{x , y \in ( I_{+} ( \gamma ) )^c }{\sum}} J_{x,y}  \big[ \tau_{\Gamma} \big( \sigma_x \big) \tau_{\Gamma} \big( \sigma_y \big) \big]  -  \underset{\Gamma \equiv \cup_i \{ \gamma^i_1 , \cdots , \gamma^i_n \} }{\underset{ y \in (I_{+} (\gamma ))^c }{\underset{x  \in  I_{+} ( \gamma )  }{\sum}}} J_{x,y}  \big[ \tau_{\Gamma} \big( \sigma_x \big) \tau_{\Gamma} \big( \sigma_y \big) \big] \\ -    \underset{\Gamma \equiv \cup_i \{ \gamma^i_1 , \cdots , \gamma^i_n \} }{\underset{x,y \in  
 \mathrm{sp} ( \gamma ) }{\sum}} J_{x,y}  \big[ \tau_{\Gamma} \big( \sigma_x \big) \tau_{\Gamma} \big( \sigma_y \big) \big]   -   \underset{\gamma \equiv \{ \gamma_1 , \cdots , \gamma_n\}}{\underset{y \in \mathrm{sp} ( \gamma )^c}{\underset{x \in \mathrm{sp} ( \gamma )}{\sum}}  }  J_{x,y}  \big[ \tau_{\Gamma} \big( \sigma_x \big) \tau_{\Gamma} \big( \sigma_y \big) \big] \\  -  \underset{\gamma \equiv \{ \gamma_1 , \cdots , \gamma_n\}}{\underset{\Gamma \equiv \cup_i \{ \gamma^i_1 , \cdots , \gamma^i_n \}}{\underset{y \in I_{-} ( \gamma ) \cup V ( \Gamma )^c}{\underset{x \in \mathrm{sp} ( \gamma )}{\sum}  }}  }J_{x,y}  \big[          \tau_{\Gamma} \big( \sigma_x \big) \tau_{\Gamma} \big( \sigma_y \big)        \big] \text{. } 
\end{align*}

\noindent From the summation above, before evaluating each instance of $\tau_{\Gamma} \big( \sigma_x \big)$ and $\tau_{\Gamma} \big( \sigma_y \big)$, observe,

\begin{align*}
    \underset{\Gamma \equiv \cup_i \{ \gamma^i_1 , \cdots , \gamma^i_n \}}{\underset{ x \in I_{-} ( \gamma ) \cup V ( \Gamma )^c }{\underset{ y \in ( I_{-} ( \gamma ) \cup V ( \Gamma )^c )^c }{\sum}}}  J_{x,y}  \big[ \tau_{\Gamma} \big( \sigma_x \big) \tau_{\Gamma} \big( \sigma_y \big) \big] =   \underset{\Gamma  \equiv \cup_i \{ \gamma^i_1 , \cdots , \gamma^i_n \}}{\underset{ x \in I_{-} ( \gamma ) \cup V ( \Gamma )^c }{\underset{ y \in I_{+} ( \gamma ) \cup \mathrm{sp} ( \gamma )  }{\sum}}}  J_{x,y}  \big[ \tau_{\Gamma} \big( \sigma_x \big) \tau_{\Gamma} \big( \sigma_y \big) \big]   \text{, }  
\end{align*}

\noindent corresponding to the summation over $y \in I_{+} \big( \gamma \big) \cup \mathrm{sp} \big( \gamma)$ and $x \in I_{-} \big( \gamma \big) \cup V \big( \Gamma \big)^c$,

\begin{align*}
  \underset{\gamma \equiv \{ \gamma_1 , \cdots , \gamma_n \} }{\underset{y \in (I_{+ } ( \gamma ) )^c}{\underset{x \in I_{+} ( \gamma )}{\sum}}} J_{xy}  \big[ \tau_{\Gamma} \big( \sigma_x \big) \tau_{\Gamma} \big( \sigma_y \big) \big]  =  \underset{\gamma \equiv \{ \gamma_1 , 
\cdots , \gamma_n\}}{\underset{\Gamma \equiv \cup_i \{ \gamma^i_1 , \cdots , \gamma^i_n \} }{\underset{y \in I_{-} ( \gamma ) \cup V ( \Gamma ) ^c}{\underset{x \in I_{+} ( \gamma )}{\sum}}} } \big[ \tau_{\Gamma} \big( \sigma_x \big) \tau_{\Gamma} \big( \sigma_y \big) \big]     \text{, }  
\end{align*}

\noindent corresponding to the summation over $x \in I_{+} \big( \gamma \big)$ and $y \in I_{-} \big( \gamma \big) \cup V \big( \Gamma \big)^c$,

\begin{align*}
        \underset{\gamma \equiv \{ \gamma_1 , \cdots , \gamma_n\} }{\underset{y \in \mathrm{sp} ( \gamma )}{\underset{x \in \mathrm{sp} ( \gamma )}{\sum}} } J_{x,y} \big[ \tau_{\Gamma} \big( \sigma_x \big) \tau_{\Gamma} \big( \sigma_y \big) ]  =      \underset{\gamma \equiv \{ \gamma_1 , \cdots , \gamma_n\} }{\underset{y \in    I_{+} ( \gamma )        }{\underset{x \in \mathrm{sp} ( \gamma )}{\sum}} }   J_{x,y} \big[ \tau_{\Gamma} \big( \sigma_x \big) \tau_{\Gamma} \big( \sigma_y \big) ] \\  +   \underset{\gamma \equiv \{ \gamma_1 , \cdots, \gamma_n \}}{\underset{\Gamma \equiv \cup_i \{ \gamma_1 , \cdots , \gamma_n\} }{\underset{y \in    I_{-} ( \gamma )  \cup V ( \Gamma )^c      }{\underset{x \in \mathrm{sp} ( \gamma )}{\sum}} }}   J_{x,y} \big[ \tau_{\Gamma} \big( \sigma_x \big) \tau_{\Gamma} \big( \sigma_y \big) ]  +  \underset{\gamma \equiv \{ \gamma_1 , \cdots , \gamma_n\} }{\underset{y \in \textbf{Z}^d}{\underset{x \in \mathrm{sp} ( \gamma )} {\sum}}} J_{x,y} \big[ \tau_{\Gamma} \big( \sigma_x \big) \tau_{\Gamma} \big( \sigma_y \big) \big]  \text{, }  
\end{align*}

\noindent corresponding to the summation over $x \in \mathrm{sp} \big( \gamma \big)$, $y \in I_{+} \big( \gamma \big)$, $y \in I_{-} \big( \gamma \big)$, and $y \in \textbf{Z}^d$. From each of the three terms in the summation above,

\begin{align*}
     \underset{\gamma \equiv \{ \gamma_1 , \cdots , \gamma_n\} }{\underset{y \in    I_{+} ( \gamma )        }{\underset{x \in \mathrm{sp} ( \gamma )}{\sum}} }   J_{x,y} \big[ \tau_{\Gamma} \big( \sigma_x \big) \tau_{\Gamma} \big( \sigma_y \big) ]    \equiv \underset{\gamma \equiv \{ \gamma_1 , \cdots , \gamma_n \}}{\underset{y \in    I_{+} ( \gamma )        }{\underset{x \in \mathrm{sp} ( \gamma )}{\sum}} }   J_{x,y}                 \Longleftrightarrow \tau_{\Gamma} \big(  \sigma_y \big) = 1 \text{, }  
\end{align*}

\noindent and $0$ otherwise, corresponding to the first term, 

\begin{align*}
    \underset{\Gamma \equiv \cup_i \{ \gamma^i_1 , \cdots , \gamma^i_n \}}{ \underset{\gamma \equiv \{ \gamma_1 , \cdots , \gamma_n\} }{\underset{y \in    I_{-} ( \gamma )        \cup V ( \Gamma )^c }{\underset{x \in \mathrm{sp} ( \gamma )}{\sum}} }}   J_{x,y} \big[ \tau_{\Gamma} \big( \sigma_x \big) \tau_{\Gamma} \big( \sigma_y \big) ]    \equiv      \underset{\Gamma \equiv \cup_i \{ \gamma^i_1 , \cdots , \gamma^i_n \}}{\underset{\gamma_1 \equiv \{ \gamma_1 , \cdots , \gamma_n\}}{\underset{ y \in I_{-} ( \gamma ) \cup V ( \Gamma )^c }{ \underset{x \in \mathrm{sp} ( \gamma )}{\sum}}  }}    J_{x,y}           \Longleftrightarrow \tau_{\Gamma} \big(  \sigma_y \big) = - 1 \text{, }  
\end{align*}

\noindent and $0$ otherwise, corresponding to the second term, and, 

\begin{align*}
     \underset{\gamma \equiv \{ \gamma_1 , \cdots , \gamma_n\} }{\underset{y \in    \textbf{Z}^d        }{\underset{x \in \mathrm{sp} ( \gamma )}{\sum}} }   J_{x,y} \big[ \tau_{\Gamma} \big( \sigma_x \big) \tau_{\Gamma} \big( \sigma_y \big) ]    \equiv \underset{\gamma \equiv  \{ \gamma_1 , \cdots , \gamma_i \}}{{\underset{y \in    \textbf{Z}^d }{\underset{x \in \mathrm{sp} ( \gamma )}{\sum}} }}   J_{x,y} \Longleftrightarrow \tau_{\Gamma} \big( \sigma_x \big) \neq \tau_{\Gamma} \big(  \sigma_y \big) 
\end{align*}

\noindent and $0$ otherwise, corresponding to the third term.

\bigskip

\noindent For the remaining terms rather than those considered above for $x \in \mathrm{sp} \big( \Gamma \big)$ and $y \in \mathrm{sp} \big( \Gamma \big)$,

\begin{align*}
     \underset{\Gamma \equiv \{ \gamma_1 , \cdots , \gamma_n\} }{\underset{y \in    \mathrm{sp} ( \gamma )        }{\underset{x \in \mathrm{sp} ( \Gamma )}{\sum}} }   J_{x,y} \big[ \tau_{\Gamma} \big( \sigma_x \big) \tau_{\Gamma} \big( \sigma_y \big) ]    \leq  \underset{\gamma \equiv  \{ \gamma_1 , \cdots , \gamma_i  \}}{{\underset{y \in    \mathrm{sp} ( \gamma )  }{\underset{x \in \mathrm{sp} ( \gamma )}{\sum}} } }  J_{x,y}           \Longleftrightarrow \tau_{\Gamma} \big( \sigma_x \big) \neq \tau_{\Gamma} \big(  \sigma_y \big) \text{, }  
\end{align*}

\noindent and $0$ otherwise. On the other hand, for the Hamiltonian of the unflipped configuration $\sigma$ that is not acted on by the mapping $\big( \tau_{\Gamma} \big( \sigma \big) \big)_x$,

\begin{align*}
    \mathcal{H}^{\mathrm{LR}, +}_{\Lambda} \big(  \sigma \big) 
 = - \underset{x,y \in \Lambda}{\sum} J_{x,y} \sigma_x \sigma_y  - \underset{y \in \Lambda^c}{\underset{x \in \Lambda}{\sum}} J_{x,y} \sigma_x \eta_y   \text{, }  
\end{align*}

\noindent from the difference,

\begin{align*}
  \mathcal{H}^{\mathrm{LR}, +}_{\Lambda} \big( \tau_{\Gamma} \big( \sigma_x \big)  \big)  -   \mathcal{H}^{\mathrm{LR}, +}_{\Lambda} \big(  \sigma \big)  =   \underset{x,y\in \Lambda}{\sum} J_{xy} \big( \tau_{\Gamma} \big( \sigma_x \big) \tau_{\Gamma} \big( \sigma_y \big) - \sigma_x \sigma_y  \big)  - \underset{y \in \Lambda^c}{\underset{x \in \Lambda}{\sum} }    J_{xy} \big( \tau_{\Gamma} \big( \sigma_x \big) \eta_y - \sigma_x \eta_y \big)    \text{, }  
\end{align*}

\noindent with $\mathcal{H}^{\mathrm{LR}, +}_{\Lambda} \big( \tau_{\Gamma} \big( \sigma \big)  \big)$ can be upper bounded with a summation over couplings, 

\begin{align*}
   \underset{y \in \mathcal{A}^{\prime}}{\underset{x \in \mathrm{sp} ( \gamma )}{\sum}} J_{x,y}  +  \underset{y \in \mathcal{B}^{\prime}}{\underset{x \in I_{-} ( \gamma )}{\sum}}  J_{x,y} +    \underset{y \in \mathcal{C}^{\prime}}{ \underset{x \in     V ( \Gamma_1 )    }{\sum} } J_{x,y}             \text{, }  
\end{align*}

\noindent which itself can be further upper bounded, as desired, by implementing the remaining argument, from \textbf{Proposition} \textit{4.5} of [2], where $\mathcal{A}^{\prime} \equiv B \big( \gamma \big)$, $\mathcal{B}^{\prime} \equiv V \big( Y_4 \big)$, $\mathcal{C}^{\prime} \equiv B \big( \gamma \big) \backslash V \big( \Gamma_2 \big)$, $\Gamma_1 \subsetneq \Gamma$, $\Gamma_2 \equiv \Gamma \backslash \Gamma_1$, and $Y_4 \equiv \Gamma_2 \backslash \big\{ \gamma^{\prime} \in \Gamma_2 :   \underset{1 \leq k \leq n}{\mathrm{sup} }\big( \mathrm{diam} \big( \gamma_k \big) \big) \leq \underset{1 \leq j \leq n }{\mathrm{sup}}   \big( \mathrm{diam} \big( \gamma^{\prime}_j \big) \big) \big\}$, in which the desired constants for the prefactor of $F_{I_{-}(\gamma)}$ are obtained from the observation that,

\begin{align*}
  \underset{y \in V ( \Gamma_{\mathrm{ext}} ( \sigma , I_{-} ( \gamma ) \backslash \{ \gamma \} )}{\underset{x \in I_{-}(\gamma)}{\sum}}  J_{x,y} \text{ } + \text{ }  \underset{y \in V ( \Gamma_{\mathrm{int}} ( \sigma , I_{-} ( \gamma )  )}{\underset{x \in I_{-}(\gamma)}{\sum}}  J_{x,y} 
 \leq         F_{I_{-}(\gamma)} \underset{> c^{\prime}_2}{\underbrace{\bigg(    \frac{2}{M^{(\alpha  - d ) \wedge 1}}    +    \frac{1}{M}    \bigg)     \kappa }} \text{, }  
\end{align*}

\noindent for realizations of exterior and interior paths, respectively given by $\Gamma_{\mathrm{ext}}$ and $\Gamma_{\mathrm{int}}$, and suitable $M,\kappa > 0$ from \textbf{Corollary} \textit{2.12} of [1], and,

\begin{align*}
  \underset{Y \in V ( \Gamma ( \sigma ) \backslash \{ \gamma \}}{ \underset{x \in \mathrm{sp} ( \gamma)}{\sum} } J_{x,y} \leq         \text{ }       \underset{> c^{\prime}_3}{\underbrace{ 2 \kappa  }}   \text{ }    F_{\mathrm{sp}(\gamma)}  \text{, }  
\end{align*}

\noindent from \textbf{Proposition} \textit{2.13} of [1], while for the remaining term, the desired upper bound takes the form,

\begin{align*}
   c^{\prime}_1 \propto           \frac{J c_{\alpha} }{\big( 2 d +1 \big) 2^{\alpha}}         \text{, }  
\end{align*}

\noindent for suitable $c_{\alpha}>0$. Hence an upper bound for the three summations above takes the form given in the proposition statement. \boxed{}

\subsection{Implementing the Ding and Zhuang approach from the upper bound in the previous section, and the coarse graining procedure}

\noindent Equipped with the upper bound of the previous section, we proceed to implement the Ding and Zhuang approach for the long range Ising model, for $d < \alpha \leq d+1$ [4], by making use of concentration results for Gaussian random variables [7]. With the results from this approach, we can upper bound the probability of bad events occurring for the long range Ising model, in the same way that bad events are upper bounded for the long range, random-field Ising model. In order to show that the probability of such bad events occurring is exponentially unlikely, we implement a three-pronged approach, consisting of steps in a Majorization of the RFIM measure, Dudley's entropy bound, and upper bounding the conditional probability under $+$ boundary conditions,

\begin{align*}
  \textbf{P}^{\mathrm{LR},+}_{\Lambda} \big[   \big|  \Delta^{\mathrm{LR}}_A \big( h \big)       \big| \geq \lambda \big| h_{A^c}      \big]  \text{, }  
\end{align*}

\noindent with a suitable exponential that is inversely proportional to the symmetric difference,

\begin{align*}
     2 \text{ } \mathrm{exp} \bigg[        - \frac{\lambda^2}{8 e^2 \big| A \big|}           \bigg]   \text{. }
\end{align*}

\bigskip

\noindent \textbf{Theorem} (\textit{it is exponentially unlikely for the complement of bad events to occur}, [6]). There exists a strictly positive constants, $C_1 \equiv C_1 \big( \alpha , d \big)$ and $\epsilon$ sufficiently large, for which,

\begin{align*}
    \textbf{P}_{\Lambda} \big[ \mathcal{B}^c \big]    \leq \mathrm{exp} \big( - C_1 \epsilon^{-2} \big)        \text{. } 
\end{align*}

\noindent \textit{Proof of Theorem}. Refer to \textbf{Proposition} \textit{3.7} of [1]. \boxed{}

\bigskip

\noindent To demonstrate that a result similar to the \textbf{Theorem} above holds, introduce similar quantities to those for the long range, random-field Ising model, namely,

\begin{align*}
   \Delta^{\mathrm{LR}}_A \big( h \big) \equiv - \frac{1}{\beta} \mathrm{log} \bigg[   \frac{Z^{+}_{\Lambda, \beta} \big( \eta \big) }{Z^{+}_{\Lambda, \beta} \big( \tau^{\mathrm{LR}}_A \big( \eta \big) \big) }     \bigg] \text{, }  
\end{align*}

\noindent for $\big( \tau^{\mathrm{LR}}_A \big( \eta \big) \big)_{\partial A} \equiv \tau^{\mathrm{LR}}_A \big( \eta \big)$, corresponding to the log-transform of the ratio of the partition functions from the long range flipping procedure applied to the boundary field $\eta$,

\begin{align*}
     \mathcal{B}^{\mathrm{LR}} \equiv   \bigg\{   \underset{\gamma \in \Gamma_0 }{\mathrm{sup}} \frac{\big|  \Delta_{I_{-} ( \gamma )} \big( \eta \big) \big|}{c^{\prime}_1  \big|  \gamma \big| }  < 1 \bigg\}      \text{, }  
\end{align*}

\noindent corresponding to the supremum of paths for which the ratio above is $< 1$, and,

\begin{align*}
       \big( \mathcal{B}^{\mathrm{LR}} \big)^c  \equiv   \bigg\{       \underset{\gamma \in \Gamma_0 }{\mathrm{sup}} \frac{\big|  \Delta_{I_{-} ( \gamma )} \big( \eta \big) \big|}{c^{\prime}_1  \big|  \gamma \big| }  > 1  \bigg\}     \text{, }  
\end{align*}

\noindent corresponding to the complement of bad events. With these quantities, to demonstrate that a result similar to the \textbf{Theorem} above holds, we make use of an entropy bound and Dudley's argument [5]. For these components of the argument, define,

\begin{align*}
  \gamma_{\theta} \big( T , d \big) \equiv \underset{( A_n )_{n \geq 0}}{\mathrm{inf}} \bigg[ \underset{t \in T}{\mathrm{sup}}  \text{ } \sum_{n \geq 0} 2^{\frac{n}{\theta}} \text{ } \mathrm{diam} \big( A_n \big( t \big) \big)   \bigg]  \text{, }  
\end{align*}

\noindent corresponding to the infimum-supremum of the summation over diameters of $A_n \big( t \big)$ for $n \geq 0$, where $A_n \big( t\big)$ denotes a partition of time, $T$, satisfying the properties:

\begin{itemize}
    \item[$\bullet$] \textit{Property 1}: The cardinality of the first partition is $\big|  A_0 \big| \equiv 1$,

    \item[$\bullet$] \textit{Property 2}: The upper bound for the cardinality of the n th partition is $\big| A_n \big| \leq 2^{2^n}$,

    \item[$\bullet$] \textit{Property 3}: The sequence of partitions $\big( A_n \big( t \big) \big)_{n \geq 0}$ is increasing, in which $A_{n+1} \big( t \big) \subsetneq A_n \big( t \big)$ for all $n$.
\end{itemize}

\noindent We will restrict our attention of the quantity above, $\gamma_{\theta} \big( T , d \big)$, for $\theta \equiv 2$.

\bigskip

\noindent In addition to these components, we implement, in order, a series of results consisting of the Majorizing measure theorem [12] (restated as \textbf{Theorem} \textit{3.9} in [1]), Dudley's entropy bound [5] (restated as \textbf{Proposition} \textit{3.10} in [1]), as well as an upper bound for the probability of the process $X_t$ obtaining a supremum exceeds a factor dependent upon $\gamma_2 \big( T , d \big)$, and on $\mathrm{diam} \big( T \big)$ [12] (restated as \textbf{Theorem} \textit{3.11} in [1]). Before implementing these three steps, we argue that a version of $\textbf{Lemma}$ \textit{1} holds for the long range Ising model, from arguments originally implemented in the case of the long range, random field Ising model.

\bigskip

\noindent \textbf{Lemma} \textit{2} (\textit{an adaptation of Lemma 1 from the Ding-Zhuang approach for the long range Ising model}, [4]). For $A , A^{\prime} \subsetneq \textbf{Z}^d$, with $A \cap A^{\prime} \neq \emptyset$ and $\big| A \big| , \big| A^{\prime} \big| < + \infty$,

\begin{align*}
  \textbf{P}^{\mathrm{LR},+}_{\Lambda} \big[   \big|  \Delta^{\mathrm{LR}}_A \big( h \big)       \big| \geq \lambda \big| h_{A^c}      \big]  \leq 2 \text{ } \mathrm{exp} \bigg[        - \frac{\lambda^2}{8 e^2 \big| A \big|}           \bigg]  \text{, }  
\end{align*}

\noindent and also that,

\begin{align*}
    \textbf{P}^{\mathrm{LR},+}_{\Lambda} \big[   \big|  \Delta^{\mathrm{LR}}_A \big( h \big) - \Delta^{\mathrm{LR}}_{A^{\prime}} \big( h \big)    \big| > \lambda \big|   h_{(A \cup A^{\prime})^c}     \big]   \leq 2 \text{ } \mathrm{exp} \bigg[     - \frac{\lambda^2}{8 e^2 \big| A \Delta A^{\prime} \big|}                \bigg]   \text{, }  
\end{align*}

\noindent for the symmetric difference between the sets $A$ and $A^{\prime}$, $A \Delta A^{\prime}$.

\bigskip

\noindent \textit{Proof of Lemma 2}. The argument directly mirrors that of \textbf{Lemma} \textit{4.1} in [4]. Initially, the primary difference arises from the fact that the $\Delta$ parameter for the long range Ising model, implying,

\begin{align*}
 \bigg|   \frac{\partial}{\partial h_{i,v}} \Delta^{\mathrm{LR}}_A \big( h \big)     \bigg| =    \bigg|       -   \frac{\sum_{\sigma}         \epsilon \sigma_v \mathrm{exp} \big( - \beta \mathcal{H}^{\mathrm{LR}} \big( \sigma \big)   \big)    }{Z^{+} \big( h \big) }      -     \frac{\sum_{\sigma}    \epsilon \sigma_v \mathrm{exp} \big( -\beta \mathcal{H}^{\mathrm{LR}} \big( \sigma \big)  
 \big)   }{Z^{+} \big( h^A \big) }        \bigg| \\   \equiv \big|   \epsilon \textbf{E}^{\mathrm{LR},+}_{\Lambda_N , \epsilon h} \big[  \sigma_v  \big] - \epsilon \textbf{E}^{\mathrm{LR},+}_{\Lambda_N , \epsilon h^A } \big[   \sigma_v  \big]   \big| \\  \equiv \big|    \epsilon    \big| \big|   \textbf{E}^{\mathrm{LR},+}_{\Lambda_N , \epsilon h} \big[  \sigma_v  \big] +  \textbf{E}^{\mathrm{LR},+}_{\Lambda_N , \epsilon h^A } \big[   \sigma_v  \big]     \big|   \\ \leq 2 \epsilon \text{, }  
\end{align*}

\noindent from which the Gaussian concentration inequality, from [7], implies the desired result for strictly positive $\epsilon$. The second inequality above can be provided with similar arguments. \boxed{}

\bigskip

\noindent Besides the result above, in order to implement the steps of the Majorizing measure theorem, Dudley's entropy bound, and an upper bound for the probability of the supremum of the process $X_t$, we provide a statement of each item used in the argument, below. 

\bigskip

\noindent \textbf{Theorem} \textit{MMT} (\textit{Majorizing measure theorem}). For a metric space $\big( T , d \big)$, and $\big( X_t \big)_{t \in T}$ with $\textbf{E} \big( X_t \big) = 0$ for every $t$, there exists some universal, strictly positive, constant $L$ for which,

\begin{align*}
  L^{-1} \gamma_2 \big( T , d \big) \leq \textbf{E} \big[ \mathrm{sup}_{t \in T} X_t  \big] \leq L \gamma_2 \big( T , d \big)   \text{. } 
\end{align*}

\noindent \textbf{Proposition} \textit{DEB} (\textit{Dudley's entropy bound}). For a family of random variables $\big( X_t \big)_{t \in T}$ satisfying,

\begin{align*}
  \textbf{P}^{\mathrm{LR},+} \big[\text{ }  \big| X_t - X_s \big| \geq \lambda  \text{ }     \big]   \leq 2 \text{ }  \mathrm{exp} \bigg(  - \big( \frac{\lambda}{\sqrt{2}} \big)^2 \big( d \big( s , t \big) \big)^{-2} \bigg) \text{, }  
\end{align*}

\noindent there exists a universal, strictly positive, constant $L$ for which,

\begin{align*}
   \textbf{E} \big[  \mathrm{sup}_{t \in T} X_t        \big] \leq L \int_0^{+\infty} \sqrt{\mathrm{log} \big[   N \big( T , d , \epsilon \big)     \big] } \text{ } \mathrm{d} \epsilon  \text{. } 
\end{align*}

\noindent \textbf{Theorem} \textit{S} (\textit{upper bounding the probability of obtaining a supremum of the process} $X_t$). For the metric space $\big( T , d \big)$, and collection $\big( X_t \big)_{t \in T}$, there exists a universal, strictly positive, constant $L$ for which,

\begin{align*}
 \textbf{P} \bigg[  \mathrm{sup}_{t \in T} X_t > L \big(  \gamma_2 \big( T , d \big) + u \text{ } \mathrm{diam} \big( T \big)       \big)      \bigg] \leq \mathrm{exp} \big( - u ^2 \big)   \text{, }  
\end{align*}

\noindent for any $u > 0$.

\bigskip

\noindent The three items above will be used to establish that the following conjecture, stated in [1], holds, which we state as another result following the next one below. 

\bigskip

\noindent

\noindent Below, we state the conjecture, and use it to prove the $\textbf{Theorem}$ for establishing that the complement of bad events occur with exponentially small probability.

\bigskip

\noindent \textbf{Conjecture} (\textit{upper bounding the probability of the complement of a bad event occurring with an exponential}, [1]). For the set of contours $\Gamma_0$ containing the origin, for any $\alpha > d$, and $d \geq 3$, there exists a constant $C_2 \equiv C_2 \big( \alpha , d \big)$ for which,

\begin{align*}
   \textbf{P} \bigg[   \underset{\gamma \in \Gamma_0}{\mathrm{sup}} \frac{\big|  \Delta_{I_{-} ( \gamma )} \big( \eta \big) \big|}{ \big|  \gamma \big| }  > 1     \bigg] \leq \mathrm{exp} \big( - C^{\prime}_2 \epsilon^{-2} \big)     \text{. } 
\end{align*}

\noindent To prove the item above, we must introduce new counting arguments for the long range contour system. To this end, we must adapt two components of the argument for proving that a phase transition occurs in the long range, random-field Ising model from [1]. Recall, from the end of \textit{2}, that the first component that the authors employ for demonstrating that the phase transition occurs is upper bounding the cardinality of,

\begin{align*}
    \mathscr{C}_l \big( \gamma \big) \equiv \underset{l \in \textbf{N}}{\bigcup} 
 \big\{ C_l : \big| C_l \cap I \big( \gamma \big) \big| \geq \frac{1}{2} \big| C_l \big| \big\}          \text{, }  
\end{align*}

\noindent which represents the set of \textit{admissible} cubes. Besides upper bounding the number of possible cubes satisfying the admissibility criteria above, the authors also upper bound the total number of paths, containing the origin and of length $n$, which is given by,

\begin{align*}
  \big| B_l \big( \Gamma_0 \big( n \big) \big)   \big| \equiv \# \big\{ \forall  C_l \text{ } ,\text{ } \exists  \gamma \in \Gamma_0 \big( n \big) :          C_l    \cap B_l    \neq \emptyset \text{ } , \text{ }      C_l \cap \gamma \neq \emptyset      \big\}   \text{, }  
\end{align*}

\bigskip

\noindent corresponding to the number of boxes covering the set of all paths containing the origin, $0$, and with length $n$. For contours that are not connected, such as those arising in long range contours, an alternative counting argument presented in [1] allows for a phase transition to be shown to occur in the long range Ising model in lower dimensions. For contours in the long range, random-field system, it was shown that an exponential upper bound on the possible number of paths can be obtained by analyzing,

\begin{align*}
  \bar{\mathscr{C}_l \big( \gamma \big) } \equiv \big\{   \forall   C_l \in \partial \mathscr{C}_l \big( \gamma \big) \text{ } \exists      \text{ }      C^{\prime}_l       :      C_l \sim C^{\prime}_l        \big\}   \text{. } 
\end{align*}

\bigskip

\noindent Below, we describe a variant of the argument provided by the authors of [1], from \textbf{Proposition} \textit{3.5}, \textbf{Proposition} \textit{3.18}, \textbf{Lemma} \textit{3.14} and \textbf{Lemma} \textit{3.17}, which we incorporate into the Dudley's entropy bound.

\bigskip

\noindent \textbf{Lemma} \textit{3} (\textit{admissibility conditions on the number of l-cubes}, \textbf{Lemma} \textit{3.14}, [1]). Fix some $A \subsetneq \textbf{Z}^d$ and $l \geq 0$. The set of admissibility criteria on the number of \textit{l-cubes}, is comprised of the two conditions,

\begin{align*}
      \frac{1}{2} \big| C_l \big| \leq \big| C_l \cap A \big|          \text{, }  \\   \big| C^{\prime}_l \cap A \big| <  \frac{1}{2} \big| C^{\prime}_l \big|     \text{, }   
\end{align*}

\noindent for the two faces $C_l$ and $C^{\prime}_l$ which overlap on exactly one face, the following lower bound holds,

\begin{align*}
2^{l(d-1)} \leq b \big|   \partial_{\mathrm{ex}} A \cap U   \big|     \text{, }  
\end{align*}

\noindent for some strictly positive $b \equiv b \big( d \big) \geq 1$.

\bigskip

\noindent In comparison to the $l$ admissiblity condition presented above from [1], a similar notion of admissiblity, $rl$ admissibility, can be used for counting the possible number of contours in the long-range RFIM. For completeness, we also provide this alternate notion of admissibility below.

\bigskip

\noindent \textbf{Lemma} \textit{4} (\textit{admissibility conditions on the number of rl-cubes}, \textbf{Lemma} \textit{3.17}, [1]). Fix some $A \subsetneq \textbf{Z}^d$, and $l \geq 0$. For the set $U \equiv C_{rl} \cup C_{r^{\prime}l}$, with $C_{rl}$ and $C_{r^{\prime}l}$ being two rl-cubes sharing exactly one face. The set of admissibility criteria is the number of $\textit{rl-cubes}$, is comprised of the two conditions,

\begin{align*}
      \frac{1}{2} \big| C_{rl} \big| \leq \big| C_{rl} \cap A \big|          \text{, }  \\   \big| C^{\prime}_{rl} \cap A \big| <  \frac{1}{2} \big| C^{\prime}_{rl} \big|     \text{, }   
\end{align*}

\noindent for the two faces $C_{rl}$ and $C^{\prime}_{rl}$ which overlap on exactly one face, the following lower bound holds,

\begin{align*}
2^{rl(d-1)} \leq b^{\prime} \big|   \partial_{\mathrm{ex}} A \cap U   \big|     \text{, }  
\end{align*}

\noindent for some strictly positive $b^{\prime} \equiv b^{\prime} \big( d \big) \geq 1$.

\bigskip

\noindent \textbf{Proposition} \textit{1} (\textit{Proposition 3.5 from} [1]). For functions the $B_0 , \cdots , B_k$, any one of which is given by,

\begin{align*}
   B_i \big( A , \textbf{Z}^d \big) \equiv B_i  \equiv   \big\{    \forall     A \subsetneq \textbf{Z}^d  \text{ } , \text{ } \exists  \text{ }  B_{\mathscr{C}_m} \equiv {\cup}_{C \in \mathscr{C}_m} C  \text{ } : \text{ }    A                 \cap C \neq \emptyset    \big\}      \text{, }  
\end{align*}

\noindent for each $1 \leq i \leq k$, there exists real constants, $b_1$ and $b_2$, with $b_1 \equiv b_1 \big( d , r \big)$, and $b_2 \equiv b_2 \big( d , r \big)$ so that,

\begin{align*}
       \big| \partial \mathscr{C}_l \big( \gamma \big) \big| \leq    b_1    \frac{\big|   \partial_{\mathrm{ex}} I \big( \gamma \big)    \big|}{2^{l(d-1)}} \leq     b_1 \frac{\big| \gamma \big|}{2^{l(d-1)}}  \text{, }  
\end{align*}

\noindent and so that,

\begin{align*}
 \big|     B_l \big( \gamma \big)      \Delta  B_{l+1} \big( \gamma \big)     \big| \leq    b_2 2^l \big| \gamma \big|       \text{. } 
\end{align*}

\bigskip

\noindent The same notions of admissibility $\textit{rl-cubes}$ can be extended to obtain an identical set of inequalities (see \textbf{Proposition} \textit{3.18} of [1]).

\bigskip

\noindent Besides the propositions above, we introduce another Proposition below for adapting Proposition \textit{3.18} from [1]. This is juxtaposed with the Entropy bound which is used to count the number of possible countours for the long rang contour system.

\bigskip

\noindent \textbf{Proposition} \textit{2} (\textit{Propoisitioon 3.18}, [1]). There exists a constant $b_4 \equiv b_4 \big( d \big)$ so that, for any natural $n$,

\begin{align*}
 \big| B_l \big( \Gamma_0 \big( n \big) \big)     \big|    \leq \mathrm{exp} \bigg(      b_4 \frac{ln}{2^{l(d-1)}}           \bigg) \text{, }  
\end{align*}

\bigskip

\noindent in which the number of coarse-grained contours contained within $B_l \big( \Gamma_0 \big( n \big) \big)$ is bounded above by an exponential.

\bigskip

\noindent For contours in the long range system, in comparison to upper bounding $B_l \big( \Gamma_0 \big( n \big) \big)$, a more complicated exponential bound, of the form stated below, also directly applies for lower dimensions of the long range Ising model. The fact that the total number of paths, of length $n$,  which contain the origin is dominated by the exponential of the length $l$ of each such path raised to some real number, in addition to the reciprocal of $2$ raised to a polynomial-logarithmic function of the dimension $d$ of the underlying lattice. For the exponential upper bound, in comparison to the notation for $B_l \big( \Gamma_0 \big( n \big) \big)$, the upper bound is for $\big| B_l \big( \mathcal{C}_0 \big( n , j \big) \big) \big|$, the number of boxes covering the set of paths,

\begin{align*}
   \mathcal{C}_0 \big( n , j \big) \equiv \big\{         \gamma \in \mathcal{E}^{+}_{\Lambda} : 0 \in V \big( \gamma \big) , \big| \gamma \big| = n     \big\}  \text{. } 
\end{align*}

\noindent \textbf{Proposition} \textit{3} (\textit{Proposition 3.31}, [1]). Fix $n,j,l \geq 0$. From the set $\mathcal{C}_0 \big( n , j\big)$ defined above, there exists a constant $c_4 \equiv c_4 \big( \alpha , d \big)$ for which,

\begin{align*}
 \big|      B_l \big( \mathcal{C}_0 \big( n , j \big) \big| \leq \mathrm{exp} \bigg(            c_4 l^k \bigg[ \frac{n}{2^{rl ( d - 1 - \frac{2\mathrm{log}_2(a) }{r-d-1-\mathrm{log}_2(a)})}} + \frac{n}{2^{2^{rl}}} + 1     \bigg]       \bigg)    \text{, }  
\end{align*}

\noindent for a suitable, strictly positive constant $a$.

\bigskip

\noindent Equipped with the counting argument for contours of the long range system, we implement the steps of the argument relying on Dudley's entropy bound, from the admissibility conditions on $\textit{rl-cubes}$.

\bigskip

\noindent \textit{Proof of Theorem and Conjecture, using Theorem S}. Applied to $\Delta_{I_{-} ( \gamma )} \big( \eta \big)$, rearranging terms after applying \textbf{Proposition} \textit{DEB} implies, for $N \equiv \mathcal{C}_0 \big(  n , j \big)$,

\begin{align*}
   \textbf{E} \bigg[   \text{ } \underset{_{\gamma \in \Gamma_0 ( n ) }   }{\mathrm{sup}}\Delta_{I_{-} ( \gamma )} \big( \eta \big) \text{ }     \bigg]    \leq  L \int_0^{+\infty} \sqrt{\mathrm{log} \big[   N \big( \mathcal{C}_0 \big( n , j \big) , d_2 , \epsilon \big)     \big] }  \text{ } \mathrm{d} \epsilon  \leq        \mathcal{C} \text{ } \overset{+\infty}{\underset{l=1}{\sum}   }            \big( 2^{\frac{rl}{2}} - 2^{\frac{rl-1}{2}} \big)     \\ \times \sqrt{\mathrm{log} \big[   N \big( \mathcal{C}_0 \big( n , j \big)  , d_2 , l^{\prime} \big)     \big] }    \text{, } 
\end{align*}

\noindent for strictly positive $\mathcal{C}$ satisfying,

\begin{align*}
   \mathcal{C}       = 2 \epsilon b_3  \sqrt{n}  \text{, } 
\end{align*}

\noindent and, for $l^{\prime} \equiv \epsilon b_3 \sqrt{2^l n}$, given in \textbf{Corollary} \textit{3.16} of [1]. From the upper bound above, we proceed to upper bound, 

\begin{align*}
\sqrt{\mathrm{log} \big[ N \big( \mathcal{C}_0 \big( n , j \big)  , d_2 , l^{\prime} \big)      \big] }   \text{, }  
\end{align*}

\noindent in which, from the counting argument for countours of the long range system that are not connected,

\begin{align*}
  \sqrt{\mathrm{log} \big[  \big| B_l \big( \mathcal{C}_0 \big( n , j \big)  \big) \big|   \big] }  \equiv \sqrt{\mathrm{log} \bigg[ \mathrm{exp} \bigg(            c_4 l^k \bigg[ \frac{n}{2^{rl ( d - 1 - \frac{2\mathrm{log}_2(a) }{r-d-1-\mathrm{log}_2(a)})}} + \frac{n}{2^{2^{rl}}} + 1     \bigg]       \bigg) \bigg] } \text{. } \end{align*}

\noindent The fact that the exponential and natural logarithm are inverse functions implies that the final expression above is equal to,
  
  \begin{align*}
  \sqrt{           c_4 l^k \bigg[ \frac{n}{2^{rl ( d - 1 - \frac{2\mathrm{log}_2(a) }{r-d-1-\mathrm{log}_2(a)})}} + \frac{n}{2^{2^{rl}}} + 1     \bigg]  }    \text{, }  
\end{align*}

\noindent hence implying,

\begin{align*}
     \mathcal{C} \text{ } \overset{+\infty}{\underset{l=1}{\sum}   }            \big( 2^{\frac{rl}{2}} - 2^{\frac{rl-1}{2}} \big)     \sqrt{\mathrm{log} \big[  \mathcal{C}_0 \big( n , j \big)  , d_2 , l^{\prime} \big)      \big] }   \leq  \mathcal{C} \overset{+\infty}{\underset{l=1}{\sum}   }            \big( 2^{\frac{rl}{2}} - 2^{\frac{rl-1}{2}} \big)          \sqrt{\mathrm{log} \big[  \big| B_l \big( \mathcal{C}_0 \big( n , j \big)  \big) \big|   \big] }   \text{, }  
  \end{align*}
  
\noindent which, in light of the previous expression obtained for $ \sqrt{\mathrm{log} \big[  \big| B_l \big( \mathcal{C}_0 \big( n , j \big)  \big) \big|   \big] }$, can be further upper bounded with,

  \begin{align*}
  \overset{+\infty}{\underset{l=1}{\sum}   }   \big( 2^{\frac{rl}{2}} - 2^{\frac{rl-1}{2}} \big)   \sqrt{           c_4 l^k \bigg[ \frac{n}{2^{rl ( d - 1 - \frac{2\mathrm{log}_2(a) }{r-d-1-\mathrm{log}_2(a)})}} + \frac{n}{2^{2^{rl}}} + 1     \bigg]  }          \text{. } 
\end{align*}

\bigskip

\noindent To remove the factors $2^{\frac{rl}{2}} - 2^{\frac{rl-1}{2}}$ for $1 \leq l \leq + \infty$ in each term of the summation in the upper bound above, observe,

\begin{align*}
  \overset{+\infty}{\underset{l=1}{\sum}   }   \big( 2^{\frac{rl}{2}} - 2^{\frac{l-1}{2}} \big) \equiv \big( \sqrt{2} - \frac{1}{\sqrt{2}} \big)  + \big(  2 - \sqrt{2}   \big)   +  \cdots  \equiv 1 - \frac{\sqrt{2}}{2}          < 1   \text{. } 
\end{align*}

\noindent This implies,

\begin{align*}
 \overset{+\infty}{\underset{l=1}{\sum}   }   \big( 2^{\frac{rl}{2}} - 2^{\frac{rl-1}{2}} \big)  \sqrt{           c_4 l^k \bigg[ \frac{n}{2^{rl ( d - 1 - \frac{2\mathrm{log}_2(a) }{r-d-1-\mathrm{log}_2(a)})}} + \frac{n}{2^{2^{rl}}} + 1     \bigg]  }       \\      \leq  \overset{+\infty}{\underset{l=1}{\sum}   } \sqrt{           c_4 l^k \bigg[ \frac{n}{2^{rl ( d - 1 - \frac{2\mathrm{log}_2(a) }{r-d-1-\mathrm{log}_2(a)})}} + \frac{n}{2^{2^{rl}}} + 1     \bigg]  }          \text{. } 
\end{align*}

\noindent Furthermore, from the upper bound above,

\begin{align*}
 \overset{+\infty}{\underset{l=1}{\sum}   } \sqrt{           c_4 l^k \bigg[ \frac{n}{2^{rl ( d - 1 - \frac{2\mathrm{log}_2(a) }{r-d-1-\mathrm{log}_2(a)})}} + \frac{n}{2^{2^{rl}}} + 1     \bigg]  }     \leq   \sqrt{c_4} \bigg[  \overset{+\infty}{\underset{l=1}{\sum}   }        \sqrt{l^k \bigg[     \frac{n}{2^{rl ( d - 1 - \frac{2\mathrm{log}_2(a) }{r-d-1-\mathrm{log}_2(a)})}} + \frac{n}{2^{2^{rl}}}    \bigg] }   \\   +    \overset{+\infty}{\underset{l=1}{\sum}   }         \sqrt{l^k}       \bigg]     \text{. } 
\end{align*}

\bigskip

\noindent From these rearrangements, one has,

\begin{align*}
  \textbf{E} \bigg[   \underset{\gamma \in \Gamma_0 ( n )}{\mathrm{sup}}    \Delta_{I_{-} ( \gamma)} \big( h \big)   \bigg] \leq \textbf{E} \bigg[    \sqrt{c_4} \bigg[  \overset{+\infty}{\underset{l=1}{\sum}   }        \sqrt{l^k \bigg[     \frac{n}{2^{rl ( d - 1 - \frac{2\mathrm{log}_2(a) }{r-d-1-\mathrm{log}_2(a)})}} + \frac{n}{2^{2^{rl}}}    \bigg] }  \\    +    \overset{+\infty}{\underset{l=1}{\sum}   }         \sqrt{l^k}       \bigg]     \bigg]    \leq b_5 \big( b_4 \big) \epsilon n            \text{. } 
\end{align*}

\noindent Before finishing the argument, first observe,

\begin{align*}
    \textbf{P} \bigg[   \underset{\gamma \in \Gamma_0}{\mathrm{sup}} \frac{\big|  \Delta_{I_{-} ( \gamma )} \big( \eta \big) \big|}{ c^{\prime}_1 \big|  \gamma \big| }  > 1       \bigg]     \approx    \textbf{P} \bigg[   \underset{\gamma \in \Gamma_0}{\mathrm{sup}} \frac{\big|  \Delta_{I_{-} ( \gamma )} \big( \eta \big) \big|}{ \big|  \gamma \big| }  > 1       \bigg]   \text{, }  
\end{align*}

\noindent from which,

\begin{align*}
  \textbf{P} \bigg[  \text{ }   \underset{\gamma \in \Gamma_0 ( n )}{\mathrm{sup}}       \frac{   \Delta_{I_{-} ( \gamma)}   \big( \eta \big)}{\big| \gamma \big| }  \geq \frac{c^{\prime}_2}{2} \text{ }   \bigg]   \equiv  \textbf{P} \bigg[  \text{ }  \underset{\gamma \in \Gamma_0 ( n )}{\mathrm{sup}}         \Delta_{I_{-} ( \gamma)}   \big( \eta \big) \geq \frac{c^{\prime}_2}{2} n     \text{ }      \bigg]   \leq \textbf{P} \bigg[     \underset{\gamma \in \Gamma_0 ( n )}{\mathrm{sup}}         \Delta_{I_{-} ( \gamma)}   \big( \eta \big) \\  \geq L \big(   b_5 \big( b_4 \big) \epsilon n  +   c^{\prime}_2   \big)        \bigg]    \text{, }  
\end{align*}

\noindent for a suitable, strictly positive, $b_5$, dependent upon $b_4$, which we achieve by applying the result,  

\begin{align*}
     \textbf{P} \bigg[  \mathrm{sup}_{t \in T} X_t > L \big(  \gamma_2 \big( T , d \big) + u \text{ } \mathrm{diam} \big( T \big)       \big)      \bigg] \leq \mathrm{exp} \big( - u ^2 \big)   \text{, }        
\end{align*}

\noindent implying the desired upper bound, upon substituting an upper bound for $\gamma_2 \big( T , d \big)$, and also for $\mathrm{diam} \big( T \big)$,

\begin{align*}
  \textbf{P} \bigg[     \underset{\gamma \in \Gamma_0 ( n )}{\mathrm{sup}}         \Delta_{I_{-} ( \gamma)}   \big( \eta \big) \geq L \big(   b_5 \big( b_4 \big) \epsilon n  - \frac{\sqrt{\mathscr{C}_2}}{\epsilon}                    \big)        \bigg] \text{, }   \end{align*}  

\noindent where,

\begin{align*}
\mathrm{diam} \big( T \big) \equiv \mathrm{diam} \big( \mathcal{C}_0 \big( n , j \big) \big)  \equiv        \underset{\gamma_1 , \gamma_2 \in \mathcal{C}_0 ( n , j )}{\mathrm{sup}}  d \big( \gamma_1 , \gamma_2 \big)   \equiv     \underset{\gamma_1 , \gamma_2 \in \mathcal{C}_0 ( n , j )}{\mathrm{sup}} \big\{  M > 0 : d \big( \gamma_1 , \gamma_2 \big) \equiv M           \big\}       \text{, }  
\end{align*}

\noindent where,

\begin{align*}
\underset{\gamma_1 , \gamma_2 \in \mathcal{C}_0 ( n , j )}{\mathrm{sup}} \big\{  M > 0 : d \big( \gamma_1 , \gamma_2 \big) \equiv M           \big\}     \propto  C \big( n , j , \epsilon , M  \big)  \big| \big| \gamma_1 - \gamma_2 \big|\big|_1   \big| I \big( \gamma_1 \big) \cap I \big( \gamma_2 \big) \big|                            \text{. }
\end{align*}

  \noindent Therefore,
  
  \begin{align*}
  \textbf{P} \bigg[ \text{ }     \underset{\gamma \in \Gamma_0 ( n )}{\mathrm{sup}}         \Delta_{I_{-} ( \gamma)}   \big( \eta \big) \geq L \bigg(   b_5 \big( b_4 \big) \epsilon n - \frac{\sqrt{\mathscr{C}_2} C^{\prime} C  }{\epsilon}       \big| \big| \gamma_1 - \gamma_2 \big|\big|_1   \big| I \big( \gamma_1 \big) \cap I \big( \gamma_2 \big) \big|      \bigg)        \bigg] \leq \mathrm{exp} \big( -   \mathscr{C}_2  \epsilon^{-2}    \big)   \text{, }  
\end{align*}

\noindent from which we conclude the argument, for suitable $\mathscr{C}_2 \equiv \mathscr{C}_2 \big( \alpha , d \big)$, and some $C \equiv C \big( n , j , \epsilon , M \big)$, $C>0$. \boxed{} 

\bigskip

\noindent We conclude with the arguments in the next section with the Peierls’ argument.

\subsection{Concluding with the classical Peierls’ argument}

\noindent In the final section, we state the inequality for executing the Peierls’ argument. The following argument is a direct adaptation of the stochastic domination provided in [1], which the probability of the event $\big\{ \sigma_0 \equiv -1 \big\}$ is dominated by the contributions of two suitable exponential functions. The first contribution of the exponential is dependent upon the inverse temperature $\beta$, while the second contribution is inversely proportional to some parameter $\epsilon$ which is taken to be sufficiently small. To obtain the desired exponential domination from two components, one argues that the intersection between the complement of bad events, $\mathcal{B}^c$, with $\big\{ \sigma_0 \equiv -1 \big\}$, can be simultaneously upper bounded together from upper bounds on each component separately. Finally, to conclude the argument, with respect to $+$ boundary conditions and the measure $\textbf{Q} \big( \cdot \big)$, the same result can be shwon to hold for the lower-dimensional RFIM measure, $\textbf{P} \big( \cdot \big)$.

\bigskip

\noindent \textbf{Theorem} (\textit{Peierls’ argument for the long range contour system, a conjecture raised in} [1]). For $d \geq 3$ and $d < \alpha \leq d+1$, there exists a suitable constant $C \equiv C \big( \alpha , d \big)$, such that, 

\begin{align*}
  \textbf{P}^{\mathrm{LR},+}_{\Lambda} \big[ \sigma_0 \equiv - 1  \big]  \leq \mathrm{exp} \big(   - C^{\prime} \beta   \big) + \mathrm{exp} \big( - C^{\prime} \epsilon^{-2} \big)  \text{, }  
\end{align*}

\noindent for the event,

\begin{align*}
  \big\{ \sigma_0 \equiv - 1 \big\}  \text{, }  
\end{align*}

\noindent for all $\beta >0$, $e \leq C^{\prime}$ and $N \geq 1$, has $\textbf{P}$-probability less than, or equal to,

\begin{align*}
 1 - \mathrm{exp} \big(   - C^{\prime} \beta   \big) + \mathrm{exp} \big( - C^{\prime} \epsilon^{-2} \big)     \text{. } 
\end{align*}

\noindent Hence, for $\beta > \beta_c$, the long range Ising model undergoes a phase transition, in which,

\begin{align*}
      \textbf{P}^{\mathrm{LR}, +}_{\Lambda , \beta, \epsilon} \neq \textbf{P}^{\mathrm{LR},-}_{\Lambda, \beta , \epsilon}     \text{, }  
\end{align*}

\noindent with $\textbf{P}$-probability $1$, as stated in $\textbf{Theorem PT}$.

\bigskip

\noindent \textit{Proof of Theorem and Theorem PT}. Under the long rang Ising model probability measure $\textbf{P}^{\mathrm{LR},+}_{\Lambda} \big( \cdot \big) \equiv \textbf{P}^{+}_{\Lambda} \big( \cdot \big)$, to demonstrate that the desired inequality holds, along the lines of the argument for \textbf{Theorem} \textit{4.1} in [1], write, from the joint probability measure,

\begin{align*}
\textbf{Q}^{\mathrm{LR},+}_{\Lambda , \beta} \big( \sigma \in A , h \in B \big)  \equiv \textbf{Q}^{+}_{\Lambda , \beta} \big( \sigma \in A , h \in B \big) \equiv \underset{B}{\int}            \textbf{P}^{\mathrm{LR},+}_{\Lambda_,\beta } \big(  A \big)   \text{ }  \mathrm{d} \textbf{P}^{\mathrm{LR},+}_{\Lambda , \beta}  \big( h \big) \\ \equiv  \underset{B}{\int}            \textbf{P}^{+}_{\Lambda_,\beta } \big(  A \big)   \text{ }  \mathrm{d} \textbf{P}^{+}_{\Lambda , \beta}  \big( h \big)   \text{, }  
\end{align*}

\noindent under $+$ boundary conditions, from which the joint probability of $\big\{ \sigma_0 \equiv -1 \big\}$,

\begin{align*}
       \textbf{Q}^{+}_{\Lambda , \beta} \big( \sigma_0 \equiv -1 \big) =    \textbf{Q}^{+}_{\Lambda , \beta} \big( \big\{ \sigma_0 \equiv -1 \big\} \cap \mathcal{B} \big)   + \textbf{Q}^{+}_{\Lambda , \beta} \big( \big\{ \sigma_0 \equiv -1 \big\} \cap \mathcal{B}^c \big)  \leq   \textbf{Q}^{+}_{\Lambda , \beta} \big( \big\{ \sigma_0 \equiv -1 \big\} \cap \mathcal{B} \big) \\ + \mathrm{exp} \big( -  C^{\prime}_1 \epsilon^{-2} \big) \text{, }  
\end{align*}

\noindent where in the last inequality, we upper bound one of the joint probability terms under $+$ boundary conditions from the fact that,

\begin{align*}
 \textbf{Q}^{+}_{\Lambda , \beta} \big( \big\{ \sigma_0 \equiv -1 \big\} \cap \mathcal{B}^c \big)  \leq   \textbf{Q}^{+}_{\Lambda , \beta} \big( \mathcal{B}^c \big)  \leq \mathrm{exp} \big( -  C^{\prime}_1 \epsilon^{-2} \big)  \text{. } 
\end{align*}

\noindent Next, write, 

\begin{align*}
 \textbf{Q}^{+}_{\Lambda , \beta} \big( \sigma_0 \equiv -1 \big)  \leq \underset{\gamma \in \mathcal{C}_0}{\sum}  \textbf{Q}^{+}_{\Lambda , \beta} \big( \Omega \big( \gamma \big)  \big)  \text{, }  
\end{align*}

\noindent corresponding to the summation over all contours $\gamma$ with $0 \in V \big( \gamma \big)$, for the collection of spins satisfying,

\begin{align*}
    \Omega \big( \gamma \big)   \equiv    \big\{ \sigma \in \Omega : \gamma \subset \Gamma \big( \sigma \big) \big\}         \text{. } 
\end{align*}

\noindent From the computations thus far with the joint measure $ \textbf{Q}^{+}_{\Lambda , \beta} \big( \cdot \big)$, we proceed to write a decomposition for,

\begin{align*}
  \textbf{Q}^{+}_{\Lambda , \beta} \big( \big\{ \sigma_0 \equiv -1 \big\} \cap \mathcal{B} \big)   \text{, }  
\end{align*}

\noindent with the integral over all possible bad events, which admits the upper bound, for, 

\begin{align*}
  \int_{\mathcal{B}}            \underset{\sigma: \sigma_0 \equiv -1}{\sum}  \mathcal{D}^{\mathrm{LR}, +}_{\Lambda , \beta} \big( \sigma , \eta \big) \mathrm{d} \eta \equiv  \int_{\mathcal{B}}            \underset{\sigma: \sigma_0 \equiv -1}{\sum}  \mathcal{D}^{+}_{\Lambda , \beta} \big( \sigma , \eta \big) \mathrm{d} \eta 
  \end{align*}
  
\noindent with, denoting $\tau^{\mathrm{LR}}_{I_{-} ( \gamma )} \big( \eta \big) \equiv \tau^{\mathrm{LR}} \big( \eta \big)$,

  \begin{align*}
  \underset{\mathcal{C}_0}{\sum}  \int_{\mathcal{B}}            \underset{ \gamma \in \sigma \in \Omega ( \gamma )}{\sum}  \mathcal{D}^{\mathrm{LR},+}_{\Lambda , \beta} \big( \sigma , \eta \big) \mathrm{d} \eta \equiv \underset{\mathcal{C}_0}{\sum}  \int_{\mathcal{B}}            \underset{ \gamma \in \sigma \in \Omega ( \gamma )}{\sum}  \mathcal{D}^{+}_{\Lambda , \beta} \big( \sigma , \eta \big) \mathrm{d} \eta \\  \leq \underset{\gamma \in \mathcal{C}_0}{\sum} \frac{2^{| \gamma |} \int_{\mathcal{B}}            \underset{ \gamma \in \sigma \in \Omega ( \gamma )}{\sum}  \mathcal{D}^{+}_{\Lambda , \beta} \big( \sigma , \eta \big) \mathrm{d} \eta   }{\int_{\mathcal{B}}            \underset{ \gamma \in \sigma \in \Omega ( \gamma )}{\sum}  \mathcal{D}^{+}_{\Lambda , \beta} \big( \tau^{\mathrm{LR}} \big( \sigma  \big) , \tau^{\mathrm{LR}} \big( \eta \big)  \big) \mathrm{d} \eta   }  \leq   \underset{\gamma \in \mathcal{C}_0}{\sum}    2^{| \gamma |} {\underset{\eta \in \mathcal{B} , \sigma \in \Omega ( \gamma ) }{\mathrm{sup}}}   \frac{\mathcal{D}^{+}_{\Lambda , \beta} \big( \sigma , \eta \big)}{\mathcal{D}^{+}_{\Lambda , \beta} \big( \tau^{\mathrm{LR}} \big( \sigma  \big) , \tau^{\mathrm{LR}} \big( \eta \big)  \big)}   \text{. } 
\end{align*}

\noindent In the rearrangements above, the $2^{| \gamma |}$ arises from the fact that,

\begin{align*}
 \int_{\mathcal{B}}         \underset{\omega \in \Omega ( \gamma )}{\sum}    \mathcal{D}^{+}_{\Lambda , \beta}   \big( \tau^{\mathrm{LR}} \big( \sigma \big) , \tau^{\mathrm{LR}} \big( \eta \big) \big) \mathrm{d} \eta       \leq  2^{| \gamma |}     \text{. } 
\end{align*}

\noindent Next, recall the identity,

\begin{align*}
\frac{\mathcal{D}^{+}_{\Lambda , \beta} \big( \sigma , \eta \big)   Z^{+ }_{\Lambda ,\beta }  \big( \eta \big) }{\mathcal{D}^{+}_{\Lambda , \beta} \big( \tau_{\gamma} \big( \sigma \big)  , \tau_{\gamma} \big( \eta \big) \big)   Z^{+}_{\Lambda ,\beta }  \big( \tau \big( \eta \big) \big)  }    =  \mathrm{exp} \big[      \beta    \mathcal{H}^{\mathrm{LR}, +}_{\Lambda} \big( \tau_{\gamma} \big( \sigma \big)  \big) - \beta  \mathcal{H}^{\mathrm{LR}, +}_{\Lambda} \big(  \sigma      \big)                 \big] \text{ }     \text{, }  
\end{align*}

\noindent and the definition of bad events $\mathcal{B}$, we proceed in the computations by upper bounding the following supremum,

\begin{align*}
     \underset{\sigma \in \Omega ( \gamma )}{ \underset{\eta \in \mathcal{B}}{ \mathrm{sup}} }           \frac{\mathcal{D}^{+}_{\Lambda , \beta} \big( \sigma , \eta \big)}{\mathcal{D}^{+}_{\Lambda , \beta} \big( \tau^{\mathrm{LR}} \big( \sigma  \big) , \tau^{\mathrm{LR}} \big( \eta \big)  \big)}    \leq \mathrm{exp} \big( - \beta c^{\prime}_2 \big| \gamma \big|  \big)    \text{ } \underset{\sigma \in \Omega ( \gamma )}{ \underset{\eta \in \mathcal{B}}{ \mathrm{sup}} }   \frac{Z^{+}_{\Lambda , \beta , \eta} \big( \tau \big( \eta \big)}{Z^{+}_{\Lambda , \beta , \eta} \big( \eta \big) }  \\  \equiv      \underset{\sigma \in \Omega ( \gamma )}{ \underset{\eta \in \mathcal{B}}{ \mathrm{sup}} }          \big[   \mathrm{exp} \big( - \beta c^{\prime}_2 \big| \gamma \big| \big) \mathrm{exp} \big( \beta \Delta_{\gamma} \big( h \big) \big) \big]         \overset{\Delta_{\gamma} ( h ) \leq  \frac{1}{2} c^{\prime}_2 | \gamma | , \forall h \in \mathcal{B} }{\leq} \mathrm{exp} \big( - \frac{\beta}{2} c^{\prime}_2 \big| \gamma \big|  \big)                    \text{. } 
\end{align*}

\noindent From the upper bound above, previous computations imply the following upper bound,

\begin{align*}
     \textbf{Q}^{+}_{\Lambda , \beta} \big( \sigma_0 \equiv -1 \big)  \leq            \underset{0 \in V ( \gamma )}{\underset{\gamma \in \mathcal{C}_0}{\sum}}     2^{| \gamma |}        \mathrm{exp} \big( - \frac{\beta}{2} c^{\prime}_2 \big| \gamma \big|  \big) + \mathrm{exp} \big( - c_0 \epsilon^{-2} \big)   \equiv    \underset{0 \in V ( \gamma )}{\underset{\gamma \in \mathcal{C}_0}{\sum}}  \mathrm{exp} \big( - \frac{\beta}{2} c^{\prime}_2 \big| \gamma \big| + \mathrm{log} 2 \big| \gamma \big| \big)  \\ + \mathrm{exp} \big( - c_0 \epsilon^{-2} \big)  \\ \leq    \underset{n \geq 1}{\underset{0 \in V ( \gamma )}{\underset{\gamma \in \mathcal{E}^{+}_{\Lambda}, | \gamma | \equiv n }{\sum}}}   \mathrm{exp} \big( - \frac{\beta}{2} c^{\prime}_2 n + \big( \mathrm{log} 2 \big) n  \big)  + \mathrm{exp} \big( - c_0 \epsilon^{-2} \big) \\ \leq \underset{n \geq 1}{\sum} \big| \mathcal{C}_0 \big( n \big) \big|  \mathrm{exp} \big( - \frac{\beta}{2} c^{\prime}_2 n + \big( \mathrm{log} 2 \big) n  \big)  + \mathrm{exp} \big( - c_0 \epsilon^{-2} \big)     \end{align*}

\noindent from which the final upper bound,

     \begin{align*}    \underset{n \geq 1}{\sum} \text{ }   \mathrm{exp}  \bigg(  \big(  C_1 - \frac{\beta}{2} c^{\prime}_2 + \mathrm{log} 2    \big) n \bigg)    + \mathrm{exp} \big( - c_0 \epsilon^{-2} \big)            \text{, }  
\end{align*}

\noindent holds, from the existence of a constant for which, 

\begin{align*}
     C_1 \geq    \frac{1}{n} \mathrm{log} \bigg[  \bigg|      \underset{n \geq 1}{\sum}   \big| \mathcal{C}_0 \big( n \big) \big|   \bigg|   \bigg]         \text{. } 
\end{align*}

\noindent Proceeding, for $\beta$ sufficiently large,

\begin{align*}
\mathrm{exp} \big( - \frac{\beta}{2} c^{\prime}_2 \big)   \leq  \mathrm{exp} \big( - 2 \beta C \big)  \text{, }  
\end{align*}

\noindent the ultimate term in the upper bound implies the following upper bound,

\begin{align*}
     \textbf{Q}^{+}_{\Lambda , \beta} \big( \sigma_0 \equiv -1 \big)   \leq   \mathrm{exp}  \big( - 2 \beta C \big)  + \mathrm{exp} \big( - c_0 \epsilon^{-2} \big)     \text{, }  
\end{align*}

\noindent for a constant satisfying,

\begin{align*}
  C \leq \frac{c^{\prime}_2}{4}  \text{. } 
\end{align*}

\noindent Altogether, we conclude the argument with the $\textbf{P}$-probability statement, in which,

\begin{align*}
  \textbf{P} \bigg[ \mathcal{D}^{+}_{\Lambda , \beta} \big( \sigma_0 \equiv -1 \big) \geq  \mathrm{exp} \big( - C \beta \big) + \mathrm{exp} \big( - C \epsilon^{-2} \big)  \bigg] \overset{(\mathrm{Markov})}{\leq}  \frac{\textbf{Q}^{+}_{\Lambda , \beta} \big( \sigma_0 \equiv -1 \big)}{\mathrm{exp} \big( - C \beta \big) + \mathrm{exp} \big( - C \epsilon^{-2} \big) } \\ \leq  \frac{\mathrm{exp} \big( - 2 \beta C \big) + \mathrm{exp} \big( - 2 C \epsilon^{-2}  \big) }{\mathrm{exp} \big( - C \beta \big) + \mathrm{exp} \big( - C \epsilon^{-2} \big) } \\ \leq \bigg(     \mathrm{exp} \big( - C \beta \big) + \mathrm{exp} \big( - C \epsilon^{-2} \big)      \bigg)^{-1}   \text{. } 
\end{align*}

\noindent Hence the desired phase transition holds with $\textbf{P}$-probability $1$, from which we conclude the argument. \boxed{}

\section{Conclusion}

In this work, we implemented a variation of the Peierls' argument, which refers to a collection of seminal arguments originally introduced by Peierls for studying phase transitions. In the lower, and higher, dimensional RFIMs alike, the phase transition is characterized by the coexistence of two probability measures under $+$, and $-$, boundary conditions. Before proving the simultaneous coexistence of two such probability measures, several characteristics of contour systems, and coarse graining, were introduced. Such objects are of significance for determining how flipping, or reversing, the orientation of countably many spins impacts the free energy, and Hamiltonian, of configurations that can be sampled over finite volume. Albeit the fact that the lower and higher dimensional RFIMs depend upon a parameter $\alpha$ that is introduced for the coupling constants $J$, configurations sampled under each measure imnply the existence of discontinuous phase transitions. Abrupt changes of behavior in physical systems that are captured through phase transitions continue to remain of interest to further explore.

\section{Declarations}

\subsection{Ethics approval and consent to participate}

The author consents to participate in the peer review process.

\subsection{Consent for publication}

The author consents to submit the following work for publication.

\subsection{Availability of data and materials}

Not applicable

\subsection{Conflict of interest}

The author declares no competing interests.

\subsection{Funding}

There are no funding sources.

\nocite{*}
\bibliography{sn-bibliography}% common bib file

\end{document}